\documentclass[a4paper]{article}
\usepackage{amssymb}
\usepackage{amsmath}

\input{tcilatex}

\newcommand{\myy}[2]{s {#1,#2}}

\begin{document}

\title{The Chow rings of generalized Grassmannians}
\author{Haibao Duan\thanks{%
Supported by NSFC.} and \ Xuezhi Zhao \\
Institute of Mathematics, Chinese Academy of Sciences, \\
dhb@math.ac.cn\\
Department of Mathematics, Capital Normal University\\
zhaoxve@mail.cnu.edu.cn\ \ \ \ \ \ \ }
\date{}
\maketitle

\begin{abstract}
Based on the basis theorem of Bruhat--Chevalley [C] and the formula for
multiplying Schubert classes obtained in [D$_{1}$] and programed in [DZ$_{{1}%
}$], we introduce a new method computing the Chow rings of flag varieties
(resp. the integral cohomology of homogeneous spaces).

The method and results of this paper have been extended in [DZ$_{3}$, DZ$%
_{4} $] to obtain the integral cohomology rings of all complete flag
manifolds, and to construct the integral cohomologies of Lie groups in terms
of Schubert classes.

\begin{description}
\item 2000 Mathematical Subject Classification: 14M15; 57T15.
\end{description}
\end{abstract}

\section{Introduction}

Let $G$ be a compact connected Lie group with a closed subgroup $H\subset G$%
. The space $G/H$ of left cosets of $H$ in $G$ is called a \textsl{%
homogeneous space}. If $H$ is the centralizer of a $1$--parameter subgroup
in $G$, $G/H$ is a smooth projective variety, called a \textsl{flag variety}.

One of the main problems in algebraic geometry (resp. topology) is to
present the Chow ring $A^{\ast }(G/H)$ of a flag variety (resp. the integral
cohomology $H^{\ast }(G/H)$ of a homogeneous space) $G/H$ by a minimal
system of generators and relations. This is a classical topic starting with
the works of H. Cartan, A. Borel, P. Baum, H. Toda and so forth. They
utilized various spectral sequence techniques for certain fibrations
associated with $G/H$ [B$_{{1}}$, B, HMS, T, Wo]. However, these techniques
encounter the same difficulties when applied to Lie groups $G$ with torsion
[I, IT, T, TW, W$_{{1}}$, W$_{{2}}$, N], in particular, when $G$ is not
prime to an exceptional Lie group.

We introduce a new method for computing the Chow ring of a flag variety
(resp. the integral cohomology of a homogeneous space). This is based on the
two results from Schubert calculus [BGG]. The first one is an additive
description of $A^{\ast }(G/H)$ in terms of Schubert classes due to
Bruhat--Chevalley [C]. The second is a formula for multiplying Schubert
classes [D$_{1}$]. Since these two results have been programed from the
Cartan matrix of $G$ in [DZ$_{{1}}$], our approach boils down the problem
directly to such primary and well known invariants of Lie groups as \textsl{%
Cartan numbers} and therefore, is self--contained in the sense that no
knowledge on the topology of Lie groups is assumed (see \S 7.1).

Starting from [DZ$_{{1}}$] our approach to the ring $A^{\ast }(G/H)$ amounts
to

\begin{quote}
a) selecting from the set of Schubert classes on $G/H$ a minimal subset that
generate $A^{\ast }(G/H)$ multiplicatively;

b) determining all non--trivial relations among these generators.
\end{quote}

\noindent We develop algebraic and computational techniques implementing
these two tasks, and demonstrate their usage in the cases of \textsl{%
generalized Grassmannians }and \textsl{rank} $1$ \textsl{homogeneous spaces}
specified below. In the subsequent works [DZ$_{3}$, DZ$_{4}$] the method and
results of this paper are extended to obtain the integral cohomology rings
of all complete flag manifolds, and to construct the integral cohomologies
of Lie groups in terms of Schubert classes.

\bigskip

Let $G$ be a Lie group with Lie algebra $L(G)$, exponential map $\exp
:L(G)\rightarrow G$, and a fixed maximal torus $T$. Let $\Omega =\{\omega _{{%
1}},\ldots ,\omega _{{n}}\}\subset L(T)$ be a set of \textsl{fundamental
dominant weights} (see\textbf{\ \S }2.1). For an $\omega \in \Omega $ the
centralizer of the $1$--parameter subgroup $H=\{\exp (t\omega )\in G\mid
t\in \mathbb{R}\}$ is called the \textsl{parabolic subgroup of }$G$\textsl{\ 
}corresponding to $\omega $. Let $H_{{s}}$ be the semi--simple part of $H$.
The flag variety $G/H$ (resp. the homogeneous space $G/H_{{s}}$) is called 
\textsl{the Grassmannian }(resp.\textsl{\ the rank} $1$ \textsl{homogeneous
space}) of $G$\ corresponding to\textsl{\ }$\omega $.

If $G$ is exceptional with rank $n$, we assume that the set $\Omega
=\{\omega _{{1}},\ldots ,\omega _{{n}}\}$ is so ordered as the
root--vertices in the Dynkin diagram of $G$ pictured in [Hu, p.58]. With
this convention we single out, for given $G$ and $\omega \in \Omega $, seven
parabolic $H$, as well as their semi--simple part $H_{{s}}$, in the table
below:

\begin{center}
\begin{tabular}{|l|l|l|l|l|l|l|l|}
\hline
$G$ & $F_{{4}}$ & $F_{{4}}$ & $E_{{6}}$ & $E_{{6}}$ & $E_{{7}}$ & $E_{{7}}$
& $E_{{8}}$ \\ \hline
$\omega $ & $\omega _{{1}}$ & $\omega _{{4}}$ & $\omega _{{2}}$ & $\omega _{{%
6}}$ & $\omega _{{1}}$ & $\omega _{{7}}$ & $\omega _{{8}}$ \\ \hline
$H$ & $C_{{3}}\cdot S^{1}$ & $B_{{3}}\cdot S^{1}$ & $A_{{6}}\cdot S^{1}$ & $%
D_{{5}}\cdot S^{1}$ & $D_{{6}}\cdot S^{1}$ & $E_{{6}}\cdot S^{1}$ & $E_{{7}%
}\cdot S^{1}$ \\ \hline
$H_{{s}}$ & $C_{{3}}$ & $B_{{3}}$ & $A_{{6}}$ & $D_{{5}}$ & $D_{{6}}$ & $E_{{%
6}}$ & $E_{{7}}$ \\ \hline
\end{tabular}
\end{center}

As applications of the methods developed in this paper, calculation is
carried out for the seven Grassmannians $G/H$ (resp. rank $1$ homogeneous
spaces $G/H_{s}$) specified above. More precisely, granted with the \textsl{%
Weyl coordinates} for Schubert classes (i.e. our indices for Schubert
classes, see \S 2.2), the following results are established.

Given a subset $\{f_{1},\ldots ,f_{m}\}$ in a ring write $\left\langle
f_{1},\ldots ,f_{m}\right\rangle $ for the ideal generated by $f_{1},\ldots
,f_{m}$.

\bigskip

\noindent \textbf{Theorem 1.} Let $y_{{1}},y_{{3}},y_{{4}},y_{{6}}$ be the
Schubert classes on $F_{{4}}/C_{{3}}\cdot S^{1}$ with Weyl coordinates $%
\sigma \lbrack 1],\sigma \lbrack 3,2,1],\sigma \lbrack 4,3,2,1],\sigma
\lbrack 3,2,4,3,2,1]$ respectively. Then

\begin{quote}
$A^{\ast }(F_{{4}}/C_{{3}}\cdot S^{1})=\mathbb{Z}[y_{{1}},y_{{3}},y_{{4}},y_{%
{6}}]/\left\langle r_{{3}},r_{{6}},r_{{8}},r_{{12}}\right\rangle $,
\end{quote}

\noindent where

\begin{quote}
$r_{{3}}=2y_{{3}}-y_{{1}}^{3}$;

$r_{{6}}=2y_{{6}}+y_{{3}}^{2}-3y_{{1}}^{2}y_{{4}}$;

$r_{{8}}=3y_{{4}}^{2}-y_{{1}}^{2}y_{{6}}$;$\quad $

$r_{{12}}=y_{{6}}^{2}-y_{{4}}^{3}$.
\end{quote}

\noindent \textbf{Theorem 2.} Let $y_{{1}},y_{{4}}$ be the Schubert classes
on $F_{{4}}/B_{{3}}\cdot S^{1}$ with Weyl coordinates $\sigma \lbrack
4],\sigma \lbrack 3,2,3,4]$ respectively, Then

\begin{center}
$A^{\ast }(F_{{4}}/B_{{3}}\cdot S^{1})=\mathbb{Z}[y_{{1}},y_{{4}%
}]/\left\langle r_{{8}},r_{{12}}\right\rangle $,
\end{center}

\noindent where

\begin{quote}
$r_{{8}}=3y_{{4}}^{2}-y_{{1}}^{8}$;

$r_{{12}}=26y_{{4}}^{3}-5y_{{1}}^{12}$.
\end{quote}

\noindent \textbf{Theorem 3.} Let $y_{{1}},y_{{3}},y_{{4}},y_{{6}}$ be the
Schubert classes on $E_{{6}}/A_{{6}}\cdot S^{1}$ with Weyl coordinates $%
\sigma \lbrack 2],\sigma \lbrack 5,4,2]\QTR{sl}{,}\sigma \lbrack {6,5,4,2}%
],\sigma \lbrack 1,3,6,5,4,2]$ respectively. Then

\begin{quote}
$A^{\ast }(E_{{6}}/A_{{6}}\cdot S^{1})=\mathbb{Z}[y_{{1}},y_{{3}},y_{{4}},y_{%
{6}}]/\left\langle r_{{6}},r_{{8}},r_{{9}},r_{{12}}\right\rangle $,
\end{quote}

\noindent where

\begin{quote}
$r_{{6}}=2y_{{6}}+y_{{3}}^{2}-3y_{{1}}^{2}y_{{4}}+2y_{{1}}^{3}y_{{3}}-y_{{1}%
}^{6}$;

$r_{{8}}=3y_{{4}}^{2}-6y_{{1}}y_{{3}}y_{{4}}+y_{{1}}^{2}y_{{6}}+5y_{{1}%
}^{2}y_{{3}}^{2}-2y_{{1}}^{5}y_{{3}}$;

$r_{{9}}=2y_{{3}}y_{{6}}-y_{{1}}^{3}y_{{6}}$;

$r_{{12}}=y_{{4}}^{3}-y_{{6}}^{2}$.
\end{quote}

\noindent \textbf{Theorem 4.} Let $y_{{1}},y_{{4}}$ be the Schubert classes
on $E_{{6}}/D_{{5}}\cdot S^{1}$ with Weyl coordinates $\sigma \lbrack
6],\sigma \lbrack 2,4,5,6]$ respectively. Then

\begin{quote}
$A^{\ast }(E_{{6}}/D_{{5}}\cdot S^{1})=\mathbb{Z}[y_{{1}},y_{{4}%
}]/\left\langle r_{{9}},r_{{12}}\right\rangle $,
\end{quote}

\noindent where

\begin{quote}
$r_{{9}}=2y_{{1}}^{9}+3y_{{1}}y_{{4}}^{2}-6y_{{1}}^{5}y_{{4}}$;

$r_{{12}}=y_{{4}}^{3}-6y_{{1}}^{4}y_{{4}}^{2}+y_{{1}}^{12}$.
\end{quote}

\noindent \textbf{Theorem 5.} Let $y_{1},y_{{5}},y_{{9}}$ be the Schubert
classes on $E_{{7}}/E_{{6}}\cdot S^{1}$ with Weyl coordinates $\sigma
\lbrack 7],\sigma \lbrack 2,4,5,6,7],\sigma \lbrack 1,5,4,2,3,4,5,6,7]$
respectively. Then

\begin{quote}
$A^{\ast }(E_{{7}}/E_{{6}}\cdot S^{1})=\mathbb{Z}[y_{{1}},y_{{5}},y_{{9}%
}]/\left\langle r_{{10}},r_{{14}},r_{{18}}\right\rangle $,
\end{quote}

\noindent where

\begin{quote}
$r_{{10}}=y_{{5}}^{2}-2y_{{1}}y_{{9}}$;

$r_{{14}}=2y_{{5}}y_{{9}}-9y_{{1}}^{4}y_{{5}}^{2}+6y_{{1}}^{9}y_{{5}}-y_{{1}%
}^{14}$;

$r_{{18}}=y_{{9}}^{2}+10y_{{1}}^{3}y_{{5}}^{3}-9y_{{1}}^{8}y_{{5}}^{2}+2y_{{1%
}}^{13}y_{{5}}$.
\end{quote}

\noindent \textbf{Theorem 6.} Let $y_{{1}},y_{{4}},y_{{6}},y_{{9}}$ be the
Schubert classes on $E_{{7}}/D_{{6}}\cdot S^{1}$ with Weyl coordinates $%
\sigma \lbrack 1],\sigma \lbrack 2,4,3,1],\sigma \lbrack 2,6,5,4,3,1],\sigma
\lbrack 3,4,2,7,6,5,4,3,1]$ respectively. Then

\begin{quote}
$A^{\ast }(E_{{7}}/D_{{6}}\cdot S^{1})=\mathbb{Z}[y_{{1}},y_{{4}},y_{{6}},y_{%
{9}}]/\left\langle r_{{9}},r_{{12}},r_{{14}},r_{{18}}\right\rangle $,
\end{quote}

\noindent where

\begin{quote}
$r_{{9}}=2y_{{9}}+3y_{{1}}y_{{4}}^{2}+4y_{{1}}^{3}y_{{6}}+2y_{{1}}^{5}y_{{4}%
}-2y_{{1}}^{9}$;

$r_{{12}}=3y_{{6}}^{2}-y_{{4}}^{3}-3y_{{1}}^{4}y_{{4}}^{2}-2y_{{1}}^{6}y_{{6}%
}+2y_{{1}}^{8}y_{{4}}$;

$r_{{14}}=3y_{{4}}^{2}y_{{6}}+3y_{{1}}^{2}y_{{6}}^{2}+6y_{{1}}^{2}y_{{4}%
}^{3}+6y_{{1}}^{4}y_{{4}}y_{{6}}+2y_{{1}}^{5}y_{{9}}-y_{{1}}^{14}$;

$r_{{18}}=5y_{{9}}^{2}+29y_{{6}}^{3}-24y_{{1}}^{6}y_{{6}}^{2}+45y_{{1}%
}^{2}y_{{4}}y_{{6}}^{2}+2y_{{1}}^{9}y_{{9}}$.
\end{quote}

\noindent \textbf{Theorem 7.} Let $y_{{1}},y_{{6}},y_{{10}},y_{{15}}$ be the
Schubert classes on $E_{{8}}/E_{{7}}\cdot S^{1}$ with Weyl coordinates $%
\sigma \lbrack 8]$, $\sigma \lbrack 3,4,5,6,7,8]$, $\sigma \lbrack
1,5,4,2,3,4,5,6,7,8]$, $\sigma \lbrack 5,4,3,1,$ $7,6,5,4,2,$ $3,4,5,6,7,8]$
respectively. Then

\begin{quote}
$A^{\ast }(E_{{8}}/E_{{7}}\cdot S^{1})=\mathbb{Z}[y_{{1}},y_{{6}},y_{{10}%
},y_{{15}}]/\left\langle r_{{15}},r_{{20}},r_{{24}},r_{{30}}\right\rangle $,
\end{quote}

\noindent where

\begin{quote}
$r_{{15}}=2y_{{15}}-16y_{{1}}^{5}y_{{10}}-10y_{{1}}^{3}y_{{6}}^{2}+10y_{{1}%
}^{9}y_{{6}}-y_{{1}}^{15}$;

$r_{{20}}=3y_{{10}}^{2}+10y_{{1}}^{2}y_{{6}}^{3}+18y_{{1}}^{4}y_{{6}}y_{{10}%
}-2y_{{1}}^{5}y_{{15}}-8y_{{1}}^{8}y_{{6}}^{2}+4y_{{1}}^{10}y_{{10}}-y_{{1}%
}^{14}y_{{6}}$;

$r_{{24}}=5y_{{6}}^{4}+30y_{{1}}^{2}y_{{6}}^{2}y_{{10}}+15y_{{1}}^{4}y_{{10}%
}^{2}-2y_{{1}}^{9}y_{{15}}-5y_{{1}}^{12}y_{{6}}^{2}+y_{{1}}^{14}y_{{10}}$;

$r_{{30}}=y_{{15}}^{2}-8y_{{10}}^{3}+y_{{6}}^{5}-2y_{{1}}^{3}y_{{6}}^{2}y_{{%
15}}+3y_{{1}}^{4}y_{{6}}y_{{10}}^{2}-8y_{{1}}^{5}y_{{10}}y_{{15}}+6y_{{1}%
}^{9}y_{{6}}y_{{15}}$

$\qquad -9y_{{1}}^{10}y_{{10}}^{2}-y_{{1}}^{12}y_{{6}}^{3}-2y_{{1}}^{14}y_{{6%
}}y_{{10}}-3y_{{1}}^{15}y_{{15}}+8y_{{1}}^{20}y_{{10}}+y_{{1}}^{24}y_{{6}%
}-y_{{1}}^{30}$.
\end{quote}

Traditionally, Schubert calculus deals with intersection theory on flag
varieties. Algorithms in \S 4.3 and the proofs of Theorems 8--14 in \S 5
demonstrate how this calculation is extended to homogeneous spaces of other
types.

This paper is so arranged. \S 2 contains a brief introduction to what we
need from Schubert calculus; \S 3 develops some algebraic results concerning
computation in the quotient of a polynomial ring. Resorting to the Gysin
sequence of circle bundles relationship between cohomologies of a
Grassmannian $G/H$ and its allied space $G/H_{{s}}$ is formulated in \S 4.
With these preliminaries, Theorems 1--7 (resp. Theorems 8--14) are
established in a unified pattern in \S 6 (resp. \S 5).

Historically, the problem of computing the Chow ring of a flag variety
(resp. the integral cohomology of a homogeneous space) has been studied by
many authors. Comparison between our method and the classical means are made
in \S 7, where mistakes occurring in the earlier computations are corrected
in \S 7.5.

Certain theoretical notion and results of this paper are also algorithmic in
nature. Their effective computability are emphasized by referring to
appropriate sections of [DZ$_{{2}}$], where intermediate data facilitating
our calculation are given in detail. To make the present work
self--contained, the most relevant data from [DZ$_{{2}}$] are summarized and
tabulated in the proofs of Theorems 8--14 in \S 5.

\section{Elements of Schubert calculus}

Assume throughout that the Lie group $G$ under consideration is compact and $%
1$--connected. Fix a maximal torus $T$ in $G$ and equip the Lie algebra $%
L(G) $ with an inner product $($ $,$ $)$, so that the adjoint representation
acts as isometries of $L(G)$. Let $\Phi =\{\beta _{{1}},\ldots ,\beta _{{n}%
}\}\subset L(T)$ be a \textsl{set of simple roots} of $G$ [Hu, p.47] (when $%
G $ is semi--simple, it is so ordered as the root--vertices in the Dynkin
diagram given in [Hu, p.58]). The \textsl{Cartan matrix} of $G$ is $C=(c_{{ij%
}})_{{n\times n}}$, where

\begin{center}
$c_{{ij}}:=2(\beta _{{i}},\beta _{{j}})/(\beta _{{j}},\beta _{{j}})$, $1\leq
i,j\leq n$ ([Hu, p.55]).
\end{center}

We recall two algorithms \textquotedblleft \textsl{Decomposition}%
\textquotedblright\ and \textquotedblleft \textsl{L--R coefficients}%
\textquotedblright\ developed in [DZ$_{{1}}$]. The first presents the Weyl
group of $G$ by the minimized decompositions of its elements, in terms of
which the Schubert varieties on $G/H$ can be constructed. The second expands
a polynomial in the Schubert classes as the linear combination of the
Schubert basis.

\textbf{2.1. Preliminaries in Weyl group. }Since $\Phi =\{\beta _{{1}%
},\ldots ,\beta _{{n}}\}$ is a basis for $L(T)$, we may introduce another
basis $\Omega =\{\omega _{{1}},\ldots ,\omega _{{n}}\}$ of $L(T)$ by

\begin{center}
$2(\omega _{{i}},\beta _{{j}})/(\beta _{{j}},\beta _{{j}})=\delta _{{i,j}}$, 
$1\leq i,j\leq n$,
\end{center}

\noindent where $\omega _{{i}}$ is known as the $i^{th}$ \textsl{fundamental
dominant weight} relative to $\Phi $ [Hu, p.67]. With respect to $\Omega $
the Cartan matrix $C=(c_{{ij}})_{{n\times n}}$ gives rise to $n$ isometries $%
\sigma _{{i}}$ on $L(T)$ by

\begin{center}
$\sigma _{{i}}(\omega _{{k}})=\left\{ 
\begin{tabular}{l}
$\omega _{{i}}\text{ if }k\neq i\text{;}$ \\ 
$\omega _{{i}}-\sum\nolimits_{{1\leq j\leq n}}c_{{ij}}\omega _{{j}}\text{ if 
}k=i$,%
\end{tabular}%
\right. $ $1\leq i\leq n$.
\end{center}

\noindent Geometrically, $\sigma _{{i}}$ is the reflection in the hyperplane 
$L_{{i}}\subset L(T)$ perpendicular to $\beta _{{i}}$ and through the origin.

Usually, the simple roots $\beta _{{i}}$'s and the fundamental weights $%
\omega _{{j}}$'s are defined as linear forms on $L(T)$ (i.e. elements in the
dual space $L(T)^{\ast }$). In this paper we identify $L(T)$ with $%
L(T)^{\ast }$ through the inner product on $L(T)$.

\bigskip

\noindent \textbf{Definition 1. }The subgroup $W(G)\subset Aut(L(T))$
generated by $\sigma _{{i}}$, $1\leq i\leq n$, is called the \textsl{Weyl
group} of $G$.

\bigskip

By Definition 1, every $w\in W(G)$ admits a factorization of the form

\begin{enumerate}
\item[(2.1)] $w=\sigma _{{i}_{{1}}}\circ \cdots \circ \sigma _{{i}_{{r}}}$, $%
1\leq i_{{1}},\ldots ,i_{{r}}\leq n$.
\end{enumerate}

\noindent Its \textsl{length} $l(w)$ is the least number of factors in all
decompositions of $w$ in the form (2.1). The decomposition (2.1) is called 
\textsl{reduced}, written $w:=\sigma \lbrack i_{{1}},\ldots ,i_{{r}}]$, if $%
r=l(w)$.

The reduced decompositions of an element $w\in W(G)$ may not be unique.
However, this ambiguity can be dispelled by employing the following notion.
Consider the set of all reduced decompositions of $w$,

\begin{center}
$D(w)=\{I=(i_{{1}},\ldots ,i_{{r}})\mid w=\sigma \lbrack I]\}$,
\end{center}

\noindent where $l(w)=r$. It can be furnished with the order $\leq $ given
by the lexicographical order on $I=(i_{{1}},\ldots ,i_{{r}})\in D(w)$. A
decomposition $w=\sigma \lbrack I]$ is called \textsl{minimized }if $I\in
D(w)$ is the minimal one with respect to $\leq $. Clearly one has

\bigskip

\noindent \textbf{Lemma 1.} \textsl{Every }$w\in W(G)$\textsl{\ has a unique
minimized decomposition.}

\bigskip

For a subset $K\subset \{1,\ldots ,n\}$ let $H_{{K}}\subset G$ be the
centralizer of the $1$--parameter subgroup $\{\exp (tb)\in G\mid t\in 
\mathbb{R}\}$, $b=\sum_{i\in K}\omega _{{i}}$. Its Weyl group $W(H_{{K}})$
is the subgroup of $W(G)$ generated by $\{\sigma _{{j}}\mid j\notin K\}$.
Resorting to the length function $l$ on $W(G)$ one may embed the set $W(H_{{K%
}};G)$ of left cosets of $W(H_{{K}})$ in $W(G)$ as the subset of $W(G)$
([BGG, 5.1]):

\begin{enumerate}
\item[(2.2)] $W(H_{{K}};G)=\{w\in W(G)\mid l(w_{{1}})\geq l(w)$, $w_{{1}}\in
wW(H_{{K}})\}$.
\end{enumerate}

\noindent We shall put $W^{r}(H_{{K}};G)=\{w\in W(H_{{K}};G)\mid l(w)=r\}$.

According to Lemma 1, every $w\in W^{r}(H_{{K}};G)$ admits a unique
minimized decomposition as $w=\sigma \lbrack I]$. As a result, with respect
to the lexicographical order on the $I$'s, $W^{r}(H_{{K}};G)$ is an ordered
set, and hence can be presented as

\begin{enumerate}
\item[(2.3)] $W^{r}(H_{{K}};G)=\{w_{{r,i}}\mid 1\leq i\leq \beta (r)\}$, $%
\beta (r):=\left\vert W^{r}(H_{{K}};G)\right\vert $,
\end{enumerate}

\noindent where $w_{{r,i}}$ is the $i^{th}$ element in $W^{r}(H_{{K}};G)$.

In [DZ$_{{1}}$] a program entitled \textquotedblleft \textsl{Decomposition}%
\textquotedblright\ has been composed, whose function is summarized below:

\begin{quote}
\textbf{Algorithm:} \textsl{Decomposition.}

\textbf{Input:} \textsl{The Cartan matrix }$C=(c_{{ij}})_{{n\times n}}$ 
\textsl{of }$G$\textsl{, and a subset }$K\subset \lbrack 1,\ldots ,n]$%
\textsl{.}

\textbf{Output: }\textsl{The set }$W(H_{{K}};G)$\textsl{\ being presented by
the minimized decompositions of its elements, together with the indexing
system (2.3) imposed by the decompositions.}
\end{quote}

\noindent \textbf{Example 1.} For those $H\subset G$ concerned by Theorems
1--7, the results from \textsl{Decomposition} are tabulated in [DZ$_{{2}}$,
1.1--7.1].$\square $

\bigskip

\textbf{2.2. Schubert varieties and the basis theorem. }Given a flag variety%
\textbf{\ }$G/H$ we can assume that the subgroup $H$ is of the form $H_{{K}}$
for some $K\subset \lbrack 1,\ldots ,n]$, since the centralizer of any $1$%
-parameter subgroup is conjugate in $G$ to one of the $H_{{K}}$ ([BH,
13.5--13.6]).

For a simple root $\beta _{{i}}\in \Phi $ let $L_{{i}}\subset L(T)$ be the
hyperplane perpendicular to $\beta _{{i}}$ and through the origin, and let $%
K_{{i}}\subset G$ be the centralizer of $\exp (L_{{i}})$. For an element $%
w\in W(H;G)$ with minimized decomposition $w=\sigma \lbrack i_{{1}},\ldots
,i_{{r}}]$, the \textsl{Schubert variety} $X_{{w}}$ on $G/H$ associated to $%
w $ is the image of the composed map

\begin{enumerate}
\item[(2.4)] $K_{{i}_{{1}}}\times \cdots \times K_{{i}_{{r}}}\rightarrow G%
\overset{p}{\rightarrow }G/H$ \ by $(k_{{1}},\ldots ,k_{{r}})\rightarrow
p(k_{{1}}\cdots k_{{r}})$,
\end{enumerate}

\noindent where $p$ is the quotient map, and where the product $\cdot $
takes place in $G$. Since

\begin{enumerate}
\item[(2.5)] \textsl{the union }$\bigcup\nolimits_{{w\in W(H;G)}}X_{{w}}$%
\textsl{\ dominates }$G/H$\textsl{\ by a cell complex with}

\textsl{\ }$\dim _{{\mathbb{R}}}X_{{w}}=2l(w)$ ([H, BGG])\textsl{,}
\end{enumerate}

\noindent one can introduce the \textsl{Schubert class} $s_{{w}}\in
A^{l(w)}(G/H)$, $w\in W(H;G)$, as the cocycle class Kronecker dual to the
fundamental cycles\textsl{\ }$[X_{{u}}]\in $ $A_{l(w)}(G/H)$ as $%
\left\langle s_{{w}},[X_{{u}}]\right\rangle =\delta _{{w,u}}$, $u\in W(H;G)$%
. (2.5) implies that ([BGG, \S 5])

\bigskip

\noindent \textbf{Lemma 2 }(\textbf{Basis theorem).} \textsl{The set of
Schubert classes }$\{s_{{w}}\mid $\textsl{\ }$w\in W(H;G)\}$\textsl{\
constitutes an additive basis for the Chow ring }$A^{\ast }(G/H)$\textsl{.}

\bigskip

Referring to the indexing system (2.3) on $W^{r}(H;G)$, we use $s_{{r,i}}$
to simplify $s_{{w}_{{r,i}}}$, and call it \textsl{the }$i^{th}$\textsl{\
Schubert class on }$G/H$\textsl{\ in degree }$r$. We create also a
definition emphasizing the role that the minimized decomposition has played
in the construction (2.4) of $X_{{w}}$ (hence of $s_{{w}}$ by duality):\ 

\bigskip

\noindent \textbf{Definition 2.} The minimized decomposition $\sigma \lbrack
I]$ of an element $w\in W(H;G)$ is called the \textsl{Weyl coordinate} of $%
s_{{w}}$.

\bigskip

\noindent \textbf{Remark 1.} a) Historically, Schubert varieties were
introduced by F. Bruhat and C. Chevalley in the context of algebraic groups
[C], whereas our definition (2.4) applies to Lie groups in the real compact
form. It was due to Hansen [H] that these two descriptions coincide.

\noindent b) In the usual definition of the Chow ring $A^{\ast }(G/H)$ the
Schubert class corresponds to a Schubert variety $X_{w}$ is a class $\alpha
_{w}$ in $A^{N-l(w)}(G/H)$, $N=\dim _{\mathbb{C}}G/H$ (instead of the
Kronecker dual $s_{{w}}$ of $[X_{w}]$ in $A^{l(w)}(G/H)$). However, via the
isomorphisms given by the cycle maps ([Fu, Example 19.1.11])

\begin{quote}
$cl:A_{r}(G/H)\cong H_{2r}(G/H)$, $cl:A^{r}(G/H)\cong H^{2r}(G/H)$
\end{quote}

\noindent and taking into account of the Poincare duality $\mathcal{D}%
:H_{2r}(G/H)\cong H^{2N-2r}(G/H)$ on $G/H$, the two elements $s_{{w}}$ and $%
\alpha _{w}$ determine each other by the relation $\alpha _{w}=\mathcal{D}%
[X_{w}]$.

\noindent c) We can use the isomorphism $cl:A^{r}(G/H)\cong H^{2r}(G/H)$ as
an identification.$\square $

\bigskip

\textbf{2.3. Multiplying Schubert classes} Let $f$ be a polynomial of degree 
$r$ in the Schubert classes $\{s_{{w}}\mid $\textsl{\ }$w\in W(H;G)\}$. By
considering $f$ as an element in $A^{r}(G/H)$ and in view of the basis
theorem, one has the expression:

\begin{enumerate}
\item[(2.6)] $f=\sum\limits_{{w\in W^{r}(H;G)}}a_{{w}}(f)s_{{w}}$, $a_{{w}%
}(f)\in \mathbb{Z}$.
\end{enumerate}

\noindent Effective computation in $A^{\ast }(G/H)$ amounts to evaluating
the integer $a_{{w}}(f)$ for every $f$ and $w$. In the special case $f=s_{{u}%
}s_{{v}}$ (i.e. a product of two Schubert classes), the integers $a_{{w}}(f)$
are well known as\textsl{\ structure constants} of $A^{\ast }(G/H)$.

Based on the multiplicative rule of Schubert classes obtained in [D$_{1}$],
a program called \textquotedblleft \textsl{Littlewood-Richardson Coefficients%
}\textquotedblright\ (abbreviated as \textsl{L-R Coefficients} in the
sequel) implementing $a_{{w}}(f)$ has been compiled in [DZ$_{{1}}$], whose
function is briefed below.

\begin{quote}
\textbf{Algorithm: }\textsl{L-R coefficients.}

\textbf{Input:} \textsl{A polynomial }$f$\textsl{\ in Schubert classes on }$%
G/H$\textsl{\ and an element }$w\in W(H;G)$\textsl{\ given by its minimized
decomposition. }

\textbf{Output: }$a_{{w}}(f)\in \mathbb{Z}$.
\end{quote}

\noindent \textbf{Example 2.} The data in [DZ$_{{2}}$, 1.2--7.2; 1.3--7.3;
1.4--7.4] are all generated by the \textsl{L-R coefficients}.$\square $

\section{The quotient of a polynomial ring}

\textbf{3.1. The problems.} In terms of the basis theorem (i.e. Lemma 2) we
may formulate our main concerns precisely.

Let $A=\oplus _{{r\geq 0}}A^{r}$ be a finitely generated graded ring. An
element $y\in A$ is called \textsl{homogeneous} of degree $r$, written $%
\left\vert y\right\vert =r$, if $y\in A^{r}$. In this paper, all elements in
a graded ring (e.g. cohomology ring; the quotient of a polynomial ring)
under consideration are homogeneous.

A subset $S=\{y_{{1}},\ldots ,y_{{n}}\}\subset A$ is called \textsl{a set of
generators} if it generates $A$ multiplicatively. A set $S$ of generators is
called \textsl{minimal }if $\left\vert S\right\vert \leq \left\vert
T\right\vert $ for any other set $T$ of generators, where $\left\vert
S\right\vert $ is the cardinality of $S$.

\begin{quote}
\textbf{Problem 1. }\textsl{For a}\textbf{\ }\textsl{given flag variety} $%
G/H $, \textsl{find a minimal set of generators of }$A^{\ast }(G/H)$\textsl{%
\ that consists of Schubert classes.}
\end{quote}

Suppose that a solution to Problem 1 is given by $S=\{y_{{1}},\ldots ,y_{{n}%
}\}$. The inclusion $S\subset A^{\ast }(G/H)$ then induces a surjective ring
map

\begin{enumerate}
\item[(3.1)] $\pi :\mathbb{Z}[y_{{1}},\ldots ,y_{{n}}]\rightarrow A^{\ast
}(G/H)$,
\end{enumerate}

\noindent whose kernel $\ker \pi \subset \mathbb{Z}[y_{{1}},\ldots ,y_{{n}}]$
is an ideal. By the Hilbert basis theorem there exists a subset $\{r_{{1}%
},\ldots ,r_{{m}}\}\subset \mathbb{Z}[y_{{1}},\ldots ,y_{{n}}]$ so that $%
\ker \pi =$\textsl{\ }$\left\langle r_{{1}},\ldots ,r_{{m}}\right\rangle $.

\begin{quote}
\textbf{Problem 2. }\textsl{Find a subset }$\{r_{{1}},\ldots ,r_{{m}%
}\}\subset \mathbb{Z}[y_{{1}},\ldots ,y_{{n}}]$\textsl{\ with }$m$\textsl{\
minimal} \textsl{so that }$\ker \pi =$\textsl{\ }$\left\langle r_{{1}%
},\ldots ,r_{{m}}\right\rangle $\textsl{.}
\end{quote}

Clearly, once both Problems 1 and 2 are solved, we arrive at the desired 
\textsl{Schubert presentation} of $A^{\ast }(G/H)$ [IM] as

\begin{enumerate}
\item[(3.2)] $A^{\ast }(G/H)=\mathbb{Z}[y_{{1}},\ldots ,y_{{n}%
}]/\left\langle r_{{1}},\ldots ,r_{{m}}\right\rangle $.
\end{enumerate}

\noindent Our investigation on Problem 1 involves geometric considerations,
and will be postponed to \S 4. This section is devoted to two algebraic
results (Lemmas 3 and 4) useful in solving Problem 2.

\bigskip

\textbf{3.2. Specifying the relations. }Let $\mathbb{Z}[y_{{1}},\ldots ,y_{{n%
}}]$ be a polynomial ring graded by $\left\vert y_{{i}}\right\vert >0$, and
let $\mathbb{Z}[y_{{1}},\ldots ,y_{{n}}]^{(m)}$ be the $\mathbb{Z}$--module
of homogeneous polynomials of degree $m$. Denote by $\mathbb{N}^{n}$ the set
of all $n$--tuples $\alpha =(b_{{1}},\ldots ,b_{{n}})$ of non--negative
integers. Then the set of monomials

\begin{enumerate}
\item[(3.3)] $B(m)=\{y^{\alpha }=y_{{1}}^{b_{1}}\cdots y_{{n}}^{b_{n}}\mid
\alpha =(b_{{1}},\ldots ,b_{{n}})\in \mathbb{N}^{n},$ $\left\vert y^{\alpha
}\right\vert =m\}$,
\end{enumerate}

\noindent forms a basis for $\mathbb{Z}[y_{{1}},\ldots ,y_{{n}}]^{(m)}$,
called the \textsl{monomial basis of }$\mathbb{Z}[y_{{1}},\ldots ,y_{{n}%
}]^{(m)}$. It will be considered as an ordered set with respect to the
lexicographical order on $\mathbb{N}^{n}$ whose cardinality is denoted by $%
b(m)$.

Given a subset $S=\{y_{{1}},\ldots ,y_{{n}}\}$ of Schubert classes on $G/H$,
let $\left\vert y_{{i}}\right\vert $ be the dimension of $y_{{i}}$ as a
cohomology class (note that the number $\left\vert y_{{i}}\right\vert $ is
always even). The inclusion $S\subset A^{\ast }(G/H)$ induces a ring map $%
\pi :\mathbb{Z}[y_{{1}},\ldots ,y_{{n}}]\rightarrow A^{\ast }(G/H)$ whose
restriction on degree $2m$ is

\begin{center}
$\pi _{{m}}:\mathbb{Z}[y_{{1}},\ldots ,y_{{n}}]^{(2m)}\rightarrow A^{m}(G/H)$%
.
\end{center}

\noindent Combining the \textsl{L--R coefficients }with the function
\textquotedblleft \textsl{Nullspace}\textquotedblright\ in \textsl{%
Mathematica}, a basis for $\ker \pi _{{m}}$ can be explicitly exhibited.

Since $A^{m}(G/H)$ has the Schubert basis $\{s_{{m,i}}\mid 1\leq i\leq $%
\textsl{\ }$\beta (m)\}$, for every $y^{\alpha }\in B(2m)$ one has a unique
expansion

\begin{center}
$\pi _{{m}}(y^{\alpha })=c_{{\alpha ,1}}s_{{m,1}}+\cdots +c_{{\alpha ,\beta
(m)}}s_{{m,\beta (m)}}$,
\end{center}

\noindent where the coefficients $c_{{\alpha ,i}}\in \mathbb{Z}$ can be
evaluated by the \textsl{L--R coefficients }as $c_{{\alpha ,i}}=a_{{w}_{{m,i}%
}}(y^{\alpha })$ (\S 2.3). The matrix $M(\pi _{{m}})=(c_{{\alpha ,i}})_{{%
b(2m)\times \beta (m)}}$ so obtained is called \textsl{the structure matrix
of }$\pi _{{m}}$.

The built--in function \textsl{Nullspace} in \textsl{Mathematica} transforms 
$M(\pi _{{m}})$ to another matrix $N(\pi _{{m}})$ in the fashion

\begin{quote}
\textsl{In:}=Nullspace[$M(\pi _{{m}})$]

\textsl{Out:}= a matrix\textsl{\ }$N(\pi _{{m}})=(b_{{j,\alpha }})_{{%
(b(2m)-\beta (m))\times b(2m)}}$
\end{quote}

\noindent whose significance is shown in the next result.

\bigskip

\noindent \textbf{Lemma 3.} \textsl{The set of polynomials }$k_{{i}%
}=\sum\limits_{y^{\alpha }\in B(2m)}b_{{i,\alpha }}y^{\alpha }$, $1\leq
i\leq b(2m)-\beta (m)$, \textsl{is a basis for }$\ker \pi _{{m}}$.$\square $

\bigskip

\noindent \textbf{Example 3. }See in [DZ$_{{2}}$, 1.4--7.4; 1.5--7.5] for
examples of structure matrices and their Nullspaces.$\square $

\bigskip

\textbf{3.3. Eliminating relations. }Let $\{r_{{1}},\ldots ,r_{{k}}\}\subset 
\mathbb{Z}[y_{{1}},\ldots ,y_{{n}}]$ be a subset. The kernel $\psi _{{m}}(r_{%
{1}},\ldots ,r_{{k}})$ of the quotient map

\begin{quote}
$\psi :\mathbb{Z}[y_{{1}},\ldots ,y_{{n}}]\rightarrow \mathbb{Z}[y_{{1}%
},\ldots ,y_{{n}}]/\left\langle r_{{1}},\ldots ,r_{{k}}\right\rangle $
\end{quote}

\noindent in degree $m$ is spanned by the subset

\begin{quote}
$\Sigma _{{m}}(r_{{1}},\ldots ,r_{{k}})=\{y^{\alpha }r_{{i}}\mid \left\vert
y^{\alpha }\right\vert +\left\vert r_{{i}}\right\vert =m\}$
\end{quote}

\noindent with cardinality

\begin{quote}
$c_{{m}}(r_{{1}},\ldots ,r_{{k}})=b(m-\left\vert r_{{1}}\right\vert )+\cdots
+b(m-\left\vert r_{{k}}\right\vert )$.
\end{quote}

\noindent On the other hand, with respect to the monomial basis $B(m)$ of%
\textsl{\ }$\mathbb{Z}[y_{{1}},\ldots ,y_{{n}}]^{(m)}$, each $y^{\alpha }r_{{%
i}}\in \Sigma _{{m}}(r_{{1}},\ldots ,r_{{k}})$ can be written uniquely as

\begin{quote}
$y^{\alpha }r_{{i}}=\underset{y^{\beta }\in B(m)}{\Sigma }a_{{(\alpha
,i),\beta }}y^{\beta }$,$\qquad a_{{(\alpha ,i),\beta }}\in \mathbb{Z}$.
\end{quote}

\noindent Write $M_{{m}}(r_{{1}},\ldots ,r_{{k}})$ for the matrix $(a_{{%
(\alpha ,i),\beta }})_{{c}_{m}{(r_{{1}},\ldots ,r_{{k}})\times b(m)}}$ (with
respect to some order on $\Sigma _{{m}}(r_{{1}},\ldots ,r_{{k}})$) so
obtained.

\bigskip

\noindent \textbf{Definition 4.} The \textsl{deficiency} $\delta _{{m}}(r_{{1%
}},\ldots ,r_{{k}})$ of the set $\{r_{{1}},\ldots ,r_{{k}}\}$ in degree $m$
is the invariant of $M_{{m}}(r_{{1}},\ldots ,r_{{k}})$ calculated by the
following procedure:

\begin{quote}
1) diagonalizing $M_{{m}}(r_{{1}},\ldots ,r_{{k}})$ using integral row and
column operations ( [S, p.163--166]);

2) setting $\delta _{{m}}(r_{{1}},\ldots ,r_{{k}})$ to be the numbers of the 
$\pm 1$'s appearing in the resulting diagonal matrix.
\end{quote}

\noindent \textbf{Example 4.} Based on the algorithm of integral row and
column reductions in [S, p.163], a program computing $\delta _{{m}}(r_{{1}%
},\ldots ,r_{{k}})$ has been composed. However, when $b(m)$ is relatively
small, $\delta _{{m}}(r_{{1}},\ldots ,r_{{k}})$ can of course be computed
directly. As an example consider in $\mathbb{Z}[y_{{1}},y_{{5}},y_{{9}}]$
with $\left\vert y_{{i}}\right\vert =2i$ the polynomials

\begin{quote}
$r_{{10}}=y_{{5}}^{2}-2y_{{1}}y_{{9}};$

$r_{{14}}=2y_{{5}}y_{{9}}-18y_{{1}}^{5}y_{{9}}+6y_{{1}}^{9}y_{{5}}-y_{{1}%
}^{14};$

$r_{{18}}=y_{{9}}^{2}+20y_{{1}}^{4}y_{{5}}y_{{9}}+2y_{{1}}^{13}y_{{5}}-18y_{{%
1}}^{9}y_{{9}}$,
\end{quote}

\noindent (see in Theorem 5). For $m=36$ we find that

\begin{quote}
$B(36)=\{{y}_{{9}}^{2},{y}_{{1}}^{3}\,{y}_{{5}}^{3},{y}_{{1}}^{4}\,y_{{5}%
}\,y_{{9}},{y}_{{1}}^{8}\,{y}_{{5}}^{2},{y}_{{1}}^{9}\,y_{{9}},{y}_{{1}%
}^{13}y_{{5}},{y}_{{1}}^{18}\}$,

$\Sigma _{{36}}(r_{{10}},r_{{14}},r_{{18}})=\{r_{{18}},y_{{1}}^{4}r_{{14}%
},y_{{1}}^{3}y_{{5}}r_{{10}},y_{{1}}^{8}r_{{10}}\}$,
\end{quote}

\noindent and that

\begin{quote}
$M_{{36}}(r_{{10}},r_{{14}},r_{{18}})=\left( 
\begin{array}{ccccccc}
1 & 0 & 20 & 0 & -18 & 2 & 0 \\ 
0 & 0 & 2 & 0 & -18 & 6 & -1 \\ 
0 & 1 & -2 & 0 & 0 & 0 & 0 \\ 
0 & 0 & 0 & 1 & -2 & 0 & 0%
\end{array}%
\right) $.
\end{quote}

\noindent These yield $b(36)=7$, $\delta _{{36}}(r_{{10}},r_{{14}},r_{{18}%
})=4$ (note that $M_{{36}}(r_{{10}},r_{{14}},r_{{18}})$ has a $4\times 4$
minor that equals to $1$).$\square $

\bigskip

For a second subset $\{g_{{1}},\ldots ,g_{{s}}\}\subset \mathbb{Z}[y_{{1}%
},\ldots ,y_{{n}}]$ we set

\begin{quote}
$A^{\ast }=\mathbb{Z}[y_{{1}},\ldots ,y_{{n}}]/\left\langle r_{{1}},\ldots
,r_{{k}},g_{{1}},\ldots ,g_{{s}}\right\rangle $
\end{quote}

\noindent and let $\varphi :\mathbb{Z}[y_{{1}},\ldots ,y_{{n}}]/\left\langle
r_{{1}},\ldots ,r_{{k}}\right\rangle \rightarrow A^{\ast }=\tbigoplus_{{%
m\geq 0}}A^{m}$ be the quotient map. The next result tells how the integers $%
b(m)$ and $\delta _{{m}}(r_{{1}},\ldots ,r_{{k}})$ are utilized in
eliminating relations on the ring $A^{\ast }$.

\bigskip

\noindent \textbf{Lemma 4. }\textsl{If rank(}$A^{m}$\textsl{)}$=b(m)-\delta
_{{m}}(r_{{1}},\ldots ,r_{{k}})$\textsl{\ for all }$m=\left\vert g_{{i}%
}\right\vert $\textsl{, }$1\leq i\leq s$\textsl{, then }$\{g_{{1}},\ldots
,g_{{s}}\}\subseteq \left\langle r_{{1}},\ldots ,r_{{k}}\right\rangle $%
\textsl{. In particular }$\varphi $\textsl{\ is a ring isomorphism.}

\noindent \textbf{Proof.} For an $1\leq i\leq s$ we set $m=\left\vert g_{{i}%
}\right\vert $, $\delta _{{m}}(r_{{1}},\ldots ,r_{{k}})=t$. Then there is a
subset $\{f_{{1}},\ldots ,f_{{t}}\}\subset \psi _{{m}}(r_{{1}},\ldots ,r_{{k}%
})$ of cardinality $t$ that can be extended to a basis $\Lambda =\{f_{{1}%
},\ldots ,f_{{t}};h_{{1}},\ldots ,h_{{b(m)-t}}\}$ of $\mathbb{Z}[y_{{1}%
},\ldots ,y_{{n}}]^{(m)}$. Consequently, one may expand $g_{{i}}$ in terms
of $\Lambda $ as

\begin{quote}
$g_{{i}}=a_{{1}}h_{{1}}+\cdots +a_{{b(m)-t}}h_{{b(m)-t}}+c_{{1}}f_{{1}%
}+\cdots +c_{{t}}f_{{t}}$, $a_{{i}},c_{{j}}\in \mathbb{Z}$.
\end{quote}

\noindent Assume on the contrary that $g_{{i}}\notin \left\langle r_{{1}%
},\ldots ,r_{{k}}\right\rangle $. Then $a_{{k}}\neq 0$ for some $1\leq k\leq
b(m)-t$. From $\varphi (g_{{i}})=0$ and $\{f_{{1}},\ldots ,f_{{t}}\}\subset $
$\left\langle r_{{1}},\ldots ,r_{{k}}\right\rangle $ we get rank($A^{m}$)$%
\leq b(m)-t-1$. This contradiction to the assumption establishes Lemma 4.$%
\square $

\bigskip

\textbf{3.4. \textsl{Giambelli polynomials} for Schubert classes.} With the
terminologies introduced in \S 3.1--3.3, it would be convenient for us to
develop here a theory of Giambelli polynomials for Schubert classes. The
result (Lemma 5) is not needed in this paper, but will play a decisive role
in extending Theorems 1--7 to flag varieties of general types ([DZ$_{{3}}$]).

For the original idea of \textsl{Giambelli polynomials }we recall the
earliest example [Hi, p.108]. If $G=U(n)$ is the unitary group of rank $n$
and if $H=U(k)\times U(n-k)$, the flag variety $G_{{n,k}}=G/H$ is the 
\textsl{Grassmannian} of $k$--planes through the origin in $\mathbb{C}^{n}$.
Let $1+c_{1}+\cdots +c_{k}$ be the total Chern class of the canonical $k$%
--bundle on $G_{{n,k}}$. Then each $c_{i}$ can be identified with an
appropriate Schubert class on $G_{{n,k}}$ (called \textsl{a} \textsl{special
Schubert class} on $G_{{n,k}}$), and with respect to these classes one has
the presentation $A^{\ast }(G_{{n,k}})=\mathbb{Z}[c_{{1}},\ldots ,c_{{k}%
}]/\left\langle r_{{n-k+1}},\ldots ,r_{{n}}\right\rangle $ in which $r_{{j}}$
is the component of the formal inverse of $1+c_{{1}}+\cdots +c_{{k}}$ in
degree $j$. It follows that every Schubert class on $G_{{n,k}}$ can be
written as a polynomial in the special ones, and such an expression is given
by the classical \textsl{Giambelli formula }[Hi, p.112].

In general, assume that $G/H$ is a flag variety and a Schubert presentation
(3.2) for $A^{\ast }(G/H)$ has been specified. As in the classical situation 
$G_{{n,k}}$ elements in $\{y_{{1}},\ldots ,y_{{n}}\}$ may be called the 
\textsl{special Schubert classes} on $G/H$, and every Schubert classes $s_{{w%
}}$ of $G/H$ can be expressed as a polynomial $\mathcal{G}_{{w}}(y_{{1}%
},\ldots ,y_{{n}})$ in these special ones. Clearly, the task of finding $%
\mathcal{G}_{{w}}$ for each $s_{{w}}$, $w\in W(H;G)$, amounts to
generalizing the classical Giambelli formula from $G_{{n,k}}$ to $G/H$.

Based on the \textsl{L--R coefficients}, a program implementing Giambelli
polynomials has been compiled, whose function is summarized below.

\begin{quote}
\textbf{Algorithm:} \textsl{Giambelli polynomials}

\textbf{Input:} \textsl{A set }$\{y_{{1}},\ldots ,y_{{n}}\}$\textsl{\ of
special Schubert classes on }$G/H$\textsl{\ and an integer }$m>0$\textsl{.}

\textbf{Output:} \textsl{Giambelli polynomials }$\mathcal{G}_{{w}}(y_{{1}%
},\ldots ,y_{{n}})$\textsl{\ for all }$w\in W^{m}(H;G)$.
\end{quote}

We clarify details in this program. With respect to the special Schubert
classes $y_{{1}},\ldots ,y_{{n}}$ write the ordered monomial basis $B(2m)$ of%
\textsl{\ }$\mathbb{Z}[y_{{1}},\ldots ,y_{{n}}]^{(2m)}$ as $%
B(2m)=\{y^{\alpha _{{1}}},\ldots ,y^{\alpha _{{b(2m)}}}\}$. The
corresponding structure matrix $M(\pi _{{m}})$ in degree $2m$ (\S 3.2) then
satisfies

\begin{center}
$\left( 
\begin{array}{c}
y^{\alpha _{{1}}} \\ 
\vdots \\ 
y^{\alpha _{{b(2m)}}}%
\end{array}%
\right) =M(\pi _{{m}})\left( 
\begin{array}{c}
s_{{m,1}} \\ 
\vdots \\ 
s_{{m,\beta (m)}}%
\end{array}%
\right) $.
\end{center}

\noindent Since $\pi _{{m}}$ is surjective, $M(\pi _{{m}})$ has a $\beta
(m)\times \beta (m)$ minor equal to $\pm 1$. The standard integral row and
column operation diagonalizing $M(\pi _{{m}})$ ([S, p.162-164]) then yields
uniquely two invertible matrices $P=P_{{b(2m)\times b(2m)}}$, $Q=Q_{{\beta
(m)\times \beta (m)}}$ that satisfy

\begin{enumerate}
\item[(3.4)] $PM(\pi _{{m}})Q=\left( 
\begin{array}{c}
I_{{\beta (m)}} \\ 
C%
\end{array}%
\right) _{{b(2m)\times \beta (m)},}$
\end{enumerate}

\noindent where $I_{{\beta (m)}}$ is the identity matrix of rank $\beta (m)$%
. The \textsl{Giambelli polynomials} is realized by the procedure below.

\begin{quote}
\textbf{Step 1.} Compute $M(\pi _{{m}})$ using the \textsl{L--R coefficients}%
;

\textbf{Step 2.} Diagonalize $M(\pi _{{m}})$ to get $P$, $Q$ in (3.4);

\textbf{Step 3.} Set $\left( 
\begin{array}{c}
\mathcal{G}_{{m,1}} \\ 
\vdots \\ 
\mathcal{G}_{{m,\beta (m)}}%
\end{array}%
\right) =Q\cdot \lbrack P^{-1}]\left( 
\begin{array}{c}
y^{\alpha _{{1}}} \\ 
\vdots \\ 
y^{\alpha _{{b(2m)}}}%
\end{array}%
\right) $,
\end{quote}

\noindent where $[P^{-1}]$ is formed by the first $\beta (m)$ rows in the
inverse $P^{-1}$ of $P$. Obviously, the polynomials $\mathcal{G}_{{m,k}}$
obtained in Step 3 depend only on the set $\{y_{{1}},\ldots ,$ $y_{{n}}\}$
of special Schubert classes. It is also clear that

\bigskip

\noindent \textbf{Lemma 5.} The $\mathcal{G}_{m,k}(y_{1},\ldots ,y_{n})$%
\textsl{\ is a Giambelli polynomial for }$s_{m,k}$\textsl{.}$\square $

\bigskip

\noindent \textbf{Example 5.} In the proof of Lemma 3 in [DZ$_{3}$] the 
\textsl{Giambelli polynomials }is applied in determining the Chern classes
of certain complex $n$-bundle on $E_{n}/A_{n}\cdot S^{1}$, $n=6,7,8$.$%
\square $

\section{\textbf{Computing with Gysin sequence}}

In this paper all cohomologies are over integer coefficients. For a
topological space $X$ we put $H^{\text{even}}(X)=\oplus _{{r\geq 0}%
}H^{2r}(X) $;\quad\ $H^{\text{odd}}(X)=\oplus _{{r\geq 0}}H^{2r+1}(X)$. Note
that $H^{\text{even}}(X)\subset H^{\ast }(X)$ is a subring.

Assume in this section that $G/H$ is the Grassmannian\textsl{\ }%
corresponding to the $k^{th}$ weight $\omega _{{k}}\in \Omega $, $1\leq
k\leq n$. From the Gysin sequence of oriented circle bundles we derive
partial solutions to problems 1 and 2 from information on $H^{\ast }(G/H_{{s}%
})$ in Lemmas 7, 8, and develop a procedure to compute $H^{\ast }(G/H_{{s}})$%
.

\textbf{4.1. Generators of the ring }$H^{\ast }(G/H)$.\textbf{\ }Since $%
W^{1}(H;G)=\{w_{{1,1}}=\sigma \lbrack k]\}$ consists of a single element,
the basis theorem implies that:

\bigskip

\noindent \textbf{Lemma 6.} $H^{2}(G/H)\cong \mathbb{Z}$\textsl{\ is
generated by }$\omega :=s_{{1,1}}$.$\square $

\bigskip

The natural projection $p:G/H_{{s}}\rightarrow G/H$ is an oriented circle
bundle on $G/H$ with Euler class $\omega $. Since $H^{\text{odd}}(G/H)=0$ by
the basis theorem, the Gysin sequence [MS, p.143] of $p$ breaks into the
short exact sequences

\begin{enumerate}
\item[(4.1)] $0\rightarrow \omega H^{2r-2}(G/H)\rightarrow H^{2r}(G/H)%
\overset{p^{\ast }}{\rightarrow }H^{2r}(G/H_{{s}})\rightarrow 0$
\end{enumerate}

\noindent as well as the isomorphisms ($\omega \cdot $ means cup--product
with $\omega $)

\begin{enumerate}
\item[(4.2)] $\beta :H^{2r-1}(G/H_{{s}})\overset{\sim }{\rightarrow }\ker
\{H^{2r-2}(G/H)\overset{\omega \cdot }{\rightarrow }H^{2r}(G/H)\}$.
\end{enumerate}

\noindent (4.1) implies that:

\textbf{\ }

\noindent \textbf{Lemma 7.} \textsl{If }$S=\{y_{{1}},\ldots ,y_{{m}%
}\}\subset H^{\ast }(G/H)$\textsl{\ is a subset so that }$p^{\ast
}S=\{p^{\ast }(y_{{1}}),\ldots ,p^{\ast }(y_{{m}})\}$\textsl{\ is a minimal
set of generators of }$H^{\text{even}}(G/H_{{s}})$\textsl{, then }$S^{\prime
}=\{\omega ,y_{{1}},\ldots ,y_{{m}}\}$\textsl{\ is a minimal set of
generators of }$H^{\ast }(G/H)$.$\square $

\bigskip

\textbf{4.2. The initial constraints of the relations.} The graded group $H^{%
\text{odd}}(G/H_{{s}})$ is always free by (4.2), and is a module over $H^{%
\text{even}}(G/H_{{s}})$ with respect to the cup--product $H^{\text{even}%
}\times H^{\text{odd}}\rightarrow H^{\text{odd}}$, $(x,y)\mapsto x\cup y$.

Let $S=\{y_{{1}},\ldots ,y_{{m}}\}$ be a subset of $H^{\ast }(G/H)$ so that $%
p^{\ast }S$ is a minimal set of generators of $H^{\text{even}}(G/H_{{s}})$.
The inclusions $\{\omega \}\cup S$ $\subset H^{\ast }(G/H)$, $p^{\ast
}S\subset H^{\ast }(G/H_{{s}})$ extend to the surjective maps $\pi $ and $%
\overline{\pi }$ that fit in the commutative diagram

\begin{enumerate}
\item[(4.3)] 
\begin{tabular}{llllll}
&  & $\mathbb{Z}[\omega ,y_{{1}},\ldots ,y_{{m}}]^{(2r)}$ & $\overset{%
\varphi }{\rightarrow }$ & $\mathbb{Z}[y_{{1}},\ldots ,y_{{m}}]^{(2r)}$ & 
\\ 
&  & $\ \ \ \pi \downarrow $ &  & $\ \ \ \overline{\pi }\downarrow $ &  \\ 
$H^{2r-2}(G/H)$ & $\overset{\omega \cdot }{\rightarrow }$ & $\ H^{2r}(G/H)$
& $\overset{p^{\ast }}{\rightarrow }$ & $\ H^{2r}(G/H_{{s}})$ & $\rightarrow
0$%
\end{tabular}
\end{enumerate}

\noindent where $\mathbb{Z}[\omega ,y_{{1}},\ldots ,y_{{m}}]$ is graded by $%
\left\vert \omega \right\vert ,\left\vert y_{{1}}\right\vert ,\ldots
,\left\vert y_{{m}}\right\vert $, and where

\begin{quote}
$\varphi (\omega )=0$, $\varphi (y_{{i}})=y_{{i}}$; $\overline{\pi }(y_{{i}%
})=p^{\ast }(y_{{i}})$.
\end{quote}

\noindent The next result allows us to formulate a partial presentation of $%
H^{\ast }(G/H)$ from information on $H^{\ast }(G/H_{{s}})$.

\bigskip

\noindent \textbf{Lemma 8. }\textsl{If }$\{h_{{1}},\ldots ,h_{{n}}\}\subset 
\mathbb{Z}[y_{{1}},\ldots ,y_{{m}}]$\textsl{\ be a subset so that}

\begin{enumerate}
\item[(4.4)] $H^{\text{even}}(G/H_{{s}})=\mathbb{Z}[p^{\ast }(y_{{1}%
}),\ldots ,p^{\ast }(y_{{m}})]/\left\langle p^{\ast }(h_{{1}}),\ldots
,p^{\ast }(h_{{n}})\right\rangle $\textsl{,}
\end{enumerate}

\noindent \textsl{and} \textsl{if }$\{d_{{1}},\ldots ,d_{{t}}\}$ \textsl{be
a basis of }$H^{\text{odd}}(G/H_{{s}})$ \textsl{over }$H^{\text{even}}(G/H_{{%
s}})$\textsl{, then for any two subsets }$\{r_{{1}},\ldots ,r_{{n}}\}$%
\textsl{, }$\{g_{{1}},\ldots ,g_{{t}}\}\subset \mathbb{Z}[\omega ,y_{{1}%
},\ldots ,y_{{m}}]$\textsl{\ that satisfy respectively \textquotedblleft the
initial constraints\textquotedblright }

\begin{quote}
1)\textsl{\ }$r_{{i}}\in \ker \pi $ \textsl{with} $r_{{i}}\mid _{{\omega =0}%
} $ $=h_{{i}}$\textsl{; }2)\textsl{\ }$\pi (g_{{j}})=\beta (d_{{j}})$\textsl{%
,}
\end{quote}

\noindent \textsl{where }$1\leq i\leq n$\textsl{,} $1\leq j\leq t$\textsl{,
one has}

\begin{enumerate}
\item[(4.5)] $H^{\ast }(G/H)=\mathbb{Z}[\omega ,y_{{1}},\ldots ,y_{{m}%
}]/\left\langle r_{{1}},\ldots ,r_{{n}};\omega g_{{1}},\ldots ,\omega g_{{t}%
}\right\rangle $\textsl{.}
\end{enumerate}

\noindent \textbf{Proof.} We begin by observing that

a) $r_{{i}}\mid _{{\omega =0}}$ $=h_{{i}}$\ is equivalent to $r_{{i}}=h_{{i}%
}+\omega f_{{i}}$\ for some $f_{{i}}$;

b) (4.4) implies that $\ker \overline{\pi }=\left\langle h_{{1}},\ldots ,h_{{%
n}}\right\rangle $.

\noindent Since $\ker \pi \supseteq \left\langle r_{{1}},\ldots ,r_{{n}%
};\omega g_{{1}},\ldots ,\omega g_{{t}}\right\rangle $, it suffices to show
that

\begin{enumerate}
\item[(4.6)] $\alpha \in \ker \pi $ implies $\alpha \in \left\langle r_{{1}%
},\ldots ,r_{{n}};\omega g_{{1}},\ldots ,\omega g_{{t}}\right\rangle $
\end{enumerate}

\noindent which will be done by induction on $2r=\left\vert \alpha
\right\vert $. The case $r=1$ is trivial by Lemma 6. Assume that (4.6) holds
for all $\alpha $ with $\left\vert \alpha \right\vert \leq 2r-2$, and
consider next the case $\left\vert \alpha \right\vert =2r$.

Write $\alpha =\alpha _{{1}}+\omega \alpha _{{2}}$, $\alpha _{{1}}\in 
\mathbb{Z}[y_{{1}},\ldots ,y_{{m}}]$. From $p^{\ast }\pi (\alpha )=\overline{%
\pi }(\alpha _{{1}})=0$ and b) we get $\alpha _{{1}}=a_{{1}}h_{{1}}+\cdots
+a_{{n}}h_{{n}}$, $a_{{i}}\in $ $\mathbb{Z}[y_{{1}},\ldots ,y_{{m}}]$. We
can rewrite in view of a) that

\begin{enumerate}
\item[(4.7)] $\alpha =a_{{1}}r_{{1}}+\cdots +a_{{n}}r_{{n}}+\omega \alpha _{{%
3}}$, where $\alpha _{{3}}=\alpha _{{2}}-(a_{{1}}f_{{1}}+\cdots +a_{{n}}f_{{n%
}})$.
\end{enumerate}

\noindent From $\pi (\alpha )=0$, $\pi (r_{{k}})=0$ we get $\pi (\alpha _{{3}%
})\in \ker \{H^{2r-2}(G/H)\overset{\cup \omega }{\rightarrow }H^{2r}(G/H)\}$%
. Since $\beta $ in (4.2) is an isomorphism, and since $H^{\text{odd}}(G/H_{{%
s}})$ is an $H^{\text{even}}(G/H_{{s}})$--module generated by $\{d_{{1}%
},\ldots ,d_{{t}}\}$ by the assumption, we obtain

\begin{enumerate}
\item[(4.8)] $\pi (\alpha _{{3}})=\tsum\limits_{1\leq i\leq t}b_{{i}}\beta
(d_{{i}})$\textsl{\ }for some\textsl{\ }$b_{{i}}\in H^{\ast }(G/H)$.
\end{enumerate}

\noindent Since $\pi $ is surjective, $b_{{i}}=\pi (q_{{i}})$ for some $q_{{i%
}}\in \mathbb{Z}[\omega ,y_{{1}},\ldots ,y_{{m}}]$. Setting $\gamma =\alpha
_{{3}}-\sum q_{{i}}g_{{i}}$, where $g_{{1}},\ldots ,g_{{t}}$ are as those in
the lemma, (4.7) becomes

\begin{enumerate}
\item[(4.9)] $\alpha =a_{{1}}r_{{1}}+\cdots +a_{{n}}r_{{n}}+\omega \gamma
+\tsum\limits_{1\leq i\leq t}b_{{i}}(\omega g_{{i}})$.
\end{enumerate}

\noindent Since $\left\vert \gamma \right\vert =\left\vert \alpha
\right\vert -2$ with $\pi (\gamma )=0$ by (4.8), the inductive hypothesis
concludes that $\gamma \in \left\langle r_{{1}},\ldots ,r_{{n}};\omega g_{{1}%
},\ldots ,\omega g_{{t}}\right\rangle $. (4.6) has now been verified by
(4.9).$\square $

\bigskip

\textbf{4.3. Algorithm for computing }$H^{\ast }(G/H_{{s}})$. We conclude
this section with a type--free procedure computing the integral cohomology
of $G/H_{{s}}$. It will be applied in the coming section to obtain $H^{\ast
}(G/H_{{s}})$ for the seven $(G,H)$ concerned by Theorems 1--7.

The procedure begins with finding an additive basis of $H^{\ast }(G/H_{{s}})$%
; followed by deriving multiplication rules for the subring $H^{\text{even}%
}(G/H_{{s}})$; and completed by describing $H^{\text{odd}}(G/H_{{s}})$ as an
module over $H^{\text{even}}(G/H_{{s}})$.

\textbf{Step 1. Finding a basis of }$H^{\ast }(G/H_{{s}})$. According to
(4.1) and (4.2) the additive groups $H^{2k-1}(G/H_{{s}})$ and $H^{2k}(G/H_{{s%
}})$ are completely determined by the homomorphism $H^{2k-2}(G/H)\overset{%
\cup \omega }{\rightarrow }H^{2k}(G/H)$. With respect to the basis $\{s_{{r,1%
}},\ldots ,$ $s_{{r,\beta (r)}}$ $\}$ of $H^{2r}(G/H)$, $\beta
(r)=\left\vert W^{r}(H;G)\right\vert $ ((2.3)), one has the expressions $%
\omega s_{{k-1,i}}=\sum a_{{i,j}}s_{{k,j}}$, $a_{{i,j}}\in \mathbb{Z}$.
Equivalently

\begin{enumerate}
\item[(4.10)] $%
\begin{array}{rcl}
\left( 
\begin{array}{l}
\omega s_{{k-1,1}} \\ 
\omega s_{{k-1,2}} \\ 
\text{ \ \ \ }\vdots \\ 
\omega s_{{k-1,\beta (k-1)}}%
\end{array}%
\right) & = & A_{{k}}\left( 
\begin{array}{l}
s_{{k,1}} \\ 
s_{{k,2}} \\ 
\text{ }\vdots \\ 
s_{{k,\beta (k)}}%
\end{array}%
\right)%
\end{array}%
$ with $A_{{k}}=(a_{{i,j}})_{{\beta (k-1)\times \beta (k)}}$.
\end{enumerate}

\noindent Since each $\omega s_{{k-1,i}}$ is a monomial in Schubert classes,
the entries of $A_{{k}}$ can be evaluated by the \textsl{L-R coefficients} (%
\S 2.3). Diagonalizing $A_{{k}}$ using the integral row and column
reductions ([S, p.162-166]) enables one to specify bases for\textsl{\ }$%
H^{2k}(G/H_{{s}})$ and $H^{2k-1}(G/H_{{s}})$\textsl{\ }in terms of Schubert
classes on $G/H$.

\bigskip

\noindent \textbf{Example 6.} For those $(G,H)$ concerned by Theorems 1--7,
the matrices $A_{{k}}$ have all been computed and tabulated in [DZ$_{{2}}$,
1.2--7.2]. See the tables in the proofs of Theorems 8--14 in \S 5 for the
basis of $H^{\ast }(G/H_{{s}})$ so derived.

\bigskip

\textbf{Step 2. The ring structure on }$H^{\text{even}}(G/H_{{s}})$. Step 1
shows how a basis of $H^{\ast }(G/H_{{s}})$ can be achieved from the
matrices $A_{{k}}$ in (4.10). In practice, in view of the surjective ring
map $p^{\ast }\colon H^{\text{even}}(G/H)\rightarrow H^{\text{even}}(G/H_{{s}%
})$, it is possible to find a subset $\Lambda $ of Schubert classes on $G/H$
so that

\begin{enumerate}
\item[(4.11)] $p^{\ast }\Lambda =\{\overline{s}_{{r,i}}:=p^{\ast }(s_{{k,i}%
})\mid s_{{k,i}}\in \Lambda \}$ is a basis of $H^{\text{even}}(G/H_{{s}})$.
\end{enumerate}

\noindent Given two basis elements $\overline{s}_{{r,i}}$, $\overline{s}_{{%
k,j}}\in p^{\ast }\Lambda $ consider their corresponding product in $H^{\ast
}(G/H)$: $s_{{r,i}}s_{{k,j}}=\sum b_{{(r,i),(k,j)}}^{t}s_{{r+k,t}}$, where
the coefficients $b_{{(r,i),(k,j)}}^{t}$ can be computed by the \textsl{L-R
coefficients}. Applying $p^{\ast }$ to this equality yields the equation in $%
H^{\text{even}}(G/H_{{s}})$

\begin{center}
$\overline{s}_{{r,i}}\overline{s}_{{k,j}}=\sum b_{{(r,i),(k,j)}}^{t}p^{\ast
}s_{{r+k,t}}$.
\end{center}

\noindent Expressing $p^{\ast }s_{{r+k,t}}$ on the right hand side in terms
of the basis elements in $p^{\ast }\Lambda $ gives rise to the \textsl{%
multiplicative rule} of $\overline{s}_{{r,i}}$ and $\overline{s}_{{k,j}}$ in 
$H^{\text{even}}(G/H_{{s}})$

\begin{enumerate}
\item[(4.12)] $\overline{s}_{{r,i}}\overline{s}_{{k,j}}=\underset{s_{{r+k,t}%
}\in \Lambda }{\sum }c_{{(r,i),(k,j)}}^{t}\overline{s}_{{r+k,t}}$
\end{enumerate}

\noindent Clearly, (4.12) suffices to characterize $H^{\text{even}}(G/H_{{s}%
})$ as a ring.

\bigskip

\noindent \textbf{Example 7.} For those $(G,H)$ concerned by Theorems 1-7,
the formulae (4.12) have all been decided and can be found in [DZ$_{{2}}$,
1.3--7.3].

\bigskip

\textbf{Step 3. }$H^{\text{odd}}(G/H_{s})$\textbf{\ as an} $H^{\text{even}%
}(G/H_{{s}})$--\textbf{module}. Since $H^{\text{odd}}$ is torsion free by
(4.2), $y\cdot H^{\text{odd}}=0$ for all $y\in $ Tor($H^{\text{even}}$). For
this reason the pairing $H^{\text{even}}\times H^{\text{odd}}\rightarrow H^{%
\text{odd}}$ in \S 4.2 reduces to a product

\begin{enumerate}
\item[(4.13)] $[H^{\text{even}}(G/H_{{s}})/$Tor($H^{\text{even}}(G/H_{{s}})$)%
$]\times H^{\text{odd}}(G/H_{{s}})\rightarrow H^{\text{odd}}(G/H_{{s}})$.
\end{enumerate}

\noindent Since $G/H_{{s}}$ is an orientable odd dimensional manifold, the
Poincar\'{e} duality implies the next result, which suffices to characterize 
$H^{\text{odd}}$ as an module over $H^{\text{even}}$ (see the proofs of
Theorems 8--14 in \S 6).

\bigskip

\noindent \textbf{Lemma 9.} \textsl{If }$\dim _{{\mathbb{R}}}G/H_{{s}}=$%
\textsl{\ }$2b+1$\textsl{, the products }(4.13)\textsl{\ in the
complementary dimensions }$[H^{2r}/$Tor($H^{2r}$)$]\times
H^{2(b-r)+1}\rightarrow H^{2b+1}=\mathbb{Z}$\textsl{\ are all non--singular.}%
$\square $

\section{Integral cohomology of $G/H_{{s}}$}

We calculate the rings $H^{\ast }(G/H_{{s}})$ for the seven homogeneous
spaces

\begin{enumerate}
\item[(5.1)] $F_{{4}}/C_{{3}}$, $F_{{4}}/B_{{3}}$, $E_{{6}}/A_{{6}}$, $E_{{6}%
}/D_{{5}}$, $E_{{7}}/D_{{6}}$, $E_{{7}}/E_{{6}}$, $E_{{8}}/E_{{7}}$.
\end{enumerate}

\noindent The results are stated in Theorems 8--14 below, where relevance of
the ring generators of $H^{\ast }(G/H_{{s}})$ with appropriate Schubert
classes on $G/H$ is emphasized.

The geometries of the spaces in (5.1) may differ considerably. However,
their cohomologies are calculated by the same procedure: following the
algorithms in \S 4.3 preliminary data produced by computer are available in
[DZ$_{{2}}$]; these are summarized in the tables contained in the proofs.
Items in these tables, together with Lemma 9, suffice to obtain $H^{\ast
}(G/H_{{s}})$.

Given a set $\{d_{{1}},\ldots ,d_{{t}}\}$ of elements graded by $\left\vert
d_{{i}}\right\vert >0$, let $\Gamma \mathbb{(}1,d_{{1}},\ldots ,d_{{t}})$ be
the graded free abelian group spanned by $1,d_{{1}},\ldots ,d_{{t}}$, and
considered as a graded ring with the trivial products $1\cdot d_{{i}}=d_{{i}%
} $; $d_{{i}}\cdot d_{{j}}=0$.

For a graded commutative ring $A$, let $A\widehat{\otimes }\Gamma (1,d_{{1}%
},\ldots ,d_{{t}})$ be the quotient of the tensor product $A\otimes \Gamma
(1,d_{{1}},\ldots ,d_{{t}})$ by the relations Tor$(A)\cdot d_{{i}}=0$, $%
1\leq i\leq t$.

We reserve the notation $s_{{r,i}}$ for the $i^{th}$ Schubert class on $G/H$
in degree $r$ (\S 2.2). If $y\in H^{\ast }(G/H)$ we write $\overline{y}%
:=p^{\ast }(y)\in H^{\ast }(G/H_{{s}})$.

\bigskip

\noindent \textbf{Theorem 8.}\ \textsl{Let }$y_{{3}}$\textsl{, }$y_{{4}}$%
\textsl{, }$y_{{6}}$\textsl{\ be the Schubert classes on }$F_{{4}}/C_{{3}%
}\cdot S^{1}$ \textsl{with Weyl coordinates }$\sigma \lbrack 3,2,1]$, $%
\sigma \lbrack 4,3,2,1]$, $\sigma \lbrack 3,2,4,3,2,1]$\textsl{\
respectively, and let }$d_{{23}}\in H^{23}(F_{{4}}/C_{{3}})$\textsl{\ be
with }$\beta (d_{{23}})=2s_{{11,1}}-s_{{11,2}}$\textsl{. Then}

\begin{center}
$H^{\ast }(F_{{4}}/C_{{3}})=\mathbb{Z}[\overline{y}_{{3}},\overline{y}_{{4}},%
\overline{y}_{{6}}]/\left\langle h_{{3}},h_{{6}},h_{{8}},h_{{12}%
}\right\rangle \widehat{\otimes }\Gamma (1,d_{{23}})$,
\end{center}

\noindent \textsl{where} $h_{{3}}=2\overline{y}_{{3}},$ $h_{{6}}=2\overline{y%
}_{{6}}+\overline{y}_{{3}}^{2}$\textsl{, }$h_{{8}}=3\overline{y}_{{4}}^{2}$%
\textsl{, }$h_{{12}}=\overline{y}_{{6}}^{2}-\overline{y}_{{4}}^{3}$\textsl{.}

\noindent \textbf{Proof. Step 1. }With the matrices $A_{{k}}$\textbf{\ }in
(4.10) being computed and presented in [DZ$_{{2}}$, 1.2], row and column
reductions yield the results in the first two columns of the table below,
which characterizes $H^{\ast }(F_{{4}}/C_{{3}})$ as a graded group:

\begin{center}
\begin{tabular}{|l|l|l|}
\hline
nontrivial $H^{k}(F_{{4}}/C_{{3}})$ & basis elements & relations \\ \hline
$H^{6}\cong \mathbb{Z}_{{2}}$ & $\bar{s}_{{3,1}}$ &  \\ \hline
$H^{8}\cong \mathbb{Z}$ & $\bar{s}_{{4,2}}$ &  \\ \hline
$H^{12}\cong \mathbb{Z}_{{4}}$ & $\bar{s}_{{6,2}}$ & $-2\bar{s}_{{6,2}}=\bar{%
s}_{{3,1}}^{2}$ \\ \hline
$H^{14}\cong \mathbb{Z}_{{2}}$ & $\bar{s}_{{7,1}}$ & $=\bar{s}_{{3,1}}%
\overline{s}_{{4,2}}$ \\ \hline
$H^{16}\cong \mathbb{Z}_{{3}}$ & $\bar{s}_{{8,1}}$ & $=-\bar{s}_{{4,2}}^{2}$
\\ \hline
$H^{18}\cong \mathbb{Z}_{{2}}$ & $\bar{s}_{{9,2}}$ & $=\bar{s}_{{3,1}}\bar{s}%
_{{6,2}}$ \\ \hline
$H^{20}\cong \mathbb{Z}_{{4}}$ & $\bar{s}_{{10,2}}$ & $=\bar{s}_{{4,2}}\bar{s%
}_{{6,2}}$ \\ \hline
$H^{26}\cong \mathbb{Z}_{{2}}$ & $\overline{s}_{{13,1}}$ & $=\bar{s}_{{3,1}}%
\bar{s}_{{4,2}}\bar{s}_{{6,2}}$ \\ \hline
$H^{23}\cong \mathbb{Z}$ & $d_{{23}}=\beta ^{-1}(2\,s_{{11,1}}-s_{{11,2}})$
&  \\ \hline
$H^{31}\cong \mathbb{Z}$ & $d_{{31}}=\beta ^{-1}(s_{{15,1}})$ & $=\pm \bar{s}%
_{{4,2}}d_{{23}}$ \\ \hline
\end{tabular}
\end{center}

\textbf{Step 2.} By the items in the second column, $H^{\text{even}}$ has a
basis of the form $p^{\ast }\Lambda $ with $\Lambda =\{s_{{3,1}},s_{{4,2}%
},s_{{6,2}},s_{{7,1}},s_{{8,1}},s_{{9,2}},s_{{10,2}},s_{{13,1}}\}$
consisting of Schubert classes on $F_{{4}}/C_{{3}}\cdot S^{1}$. Following
the algorithm given in Step 2 in \S 4.3, the multiplicative rule (4.12) for
elements in $p^{\ast }\Lambda $ have been determined ([DZ$_{{2}}$, 1.3]),
and is recorded in the last column of the table corresponding to $H^{\text{%
even}}$. These imply that, if we put $y_{{3}}=s_{{3,1}}$, $y_{{4}}=s_{{4,2}}$%
{, }$y_{{6}}=s_{{6,2}}$, then

\begin{quote}
a) $y_{{3}}$, $y_{{4}}${, }$y_{{6}}$ are the Schubert classes whose Weyl
coordinates are given as those in the theorem by [DZ$_{{2}}$, 1.1];

b) $H^{even}(F_{{4}}/C_{{3}})$ is generated by $\overline{y}_{{3}},\overline{%
y}_{{4}}{,}\overline{y}_{{6}}$ subject to $h_{{3}},h_{{6}}$ and $h_{{8}}$.
\end{quote}

\noindent Combining these with the obvious relations $\overline{y}_{{6}}^{2}=%
\overline{y}_{{4}}^{3}\in H^{24}=0$, together with the fact that $%
\left\langle h_{{3}},h,h_{{8}},\overline{y}_{{6}}^{2},\overline{y}_{{4}%
}^{3}\right\rangle =\left\langle h_{{3}},h_{{6}},h_{{8}},h_{{12}%
}\right\rangle $ in $\mathbb{Z}[\overline{y}_{{3}},\overline{y}_{{4}}{,}%
\overline{y}_{{6}}]$, one obtains

\begin{enumerate}
\item[(5.2)] $H^{\text{even}}(F_{{4}}/C_{{3}})=\mathbb{Z}[\overline{y}_{{3}},%
\overline{y}_{{4}},\overline{y}_{{6}}]/\left\langle h_{{3}},h_{{6}},h_{{8}%
},h_{{12}}\right\rangle $.
\end{enumerate}

\textbf{Step 3.} The proof is completed by $d_{{31}}=\pm \bar{s}_{{4,2}}d_{{%
23}}$ (Lemma 9) and $d_{{23}}^{2}\in H^{46}=0$ (see in the first column of
the table).$\square $

\bigskip

\noindent \textbf{Theorem 9. }\textsl{Let }$y_{{4}}$\textsl{\ be the
Schubert class on }$F_{{4}}/B_{{3}}\cdot S^{1}$ \textsl{with Weyl coordinate 
}$\sigma \lbrack 3,2,3,4]$\textsl{; and let }$d_{{23}}\in H^{23}(F_{{4}}/B_{{%
3}})$\textsl{\ be with }$\beta (d_{{23}})=-s_{{11,1}}+s_{{11,2}}$\textsl{.
Then}

\begin{center}
$H^{\ast }(F_{{4}}/B_{{3}})=\mathbb{Z}[\overline{y}_{{4}}]/\left\langle h_{{8%
}},h_{{12}}\right\rangle \widehat{\otimes }\Gamma (1,d_{{23}})$\textsl{, }
\end{center}

\noindent \textsl{where }$h_{{8}}=3\overline{y}_{{4}}^{2}$\textsl{, }$h_{{12}%
}=\overline{y}_{{4}}^{3}$\textsl{.}

\noindent \textbf{Proof. Step 1. }From the matrices\textbf{\ }$A_{{k}}$
presented in [DZ$_{{2}}$, 2.2], one deduces the results in the first two
columns of the table below.

\begin{center}
\begin{tabular}{|l|l|l|}
\hline
nontrivial $H^{k}$ & basis elements & relations \\ \hline
$H^{8}\cong \mathbb{Z}$ & $\bar{s}_{{4,2}}$ &  \\ \hline
$H^{16}\cong \mathbb{Z}_{{3}}$ & $\bar{s}_{{8,1}}$ & $={\bar{s}}_{{4,2}}^{2}$
\\ \hline
$H^{23}\cong \mathbb{Z}$ & $d_{{23}}=\beta ^{-1}(-s_{{11,1}}+s_{{11,2}})$ & 
\\ \hline
$H^{31}\cong \mathbb{Z}$ & $d_{{31}}=\beta ^{-1}s_{{15,1}}$ & $=\pm \bar{s}_{%
{4,2}}d_{{23}}$ \\ \hline
\end{tabular}
\end{center}

\textbf{Step 2.} By the results in the second column, $H^{\text{even}}(F_{{4}%
}/B_{{3}})$ has the basis $p^{\ast }\Lambda $ with $\Lambda =\{s_{{4,2}},s_{{%
8,1}}\}$. The corresponding (4.12) consists of the single equation $%
\overline{s}_{{8,1}}=\overline{{s}}_{{4,2}}^{2}$ (see [DZ$_{{2}}$, 2.3]).
These implies that, if we put $y_{{4}}=s_{{4,2}}$, then $y_{{4}}$ is the
Schubert class whose Weyl coordinate is given as those in the theorem by [DZ$%
_{{2}}$, 2.1]; and the ring $H^{\text{even}}(F_{{4}}/B_{{3}})$ is generated
by $\overline{y}_{{4}}$ subject to the relations $h_{{8}}$ and $h_{{12}}=%
\overline{y}_{{4}}^{3}$ (since $H^{24}=0$):

\begin{enumerate}
\item[(5.3)] $H^{\text{even}}(F_{{4}}/B_{{3}})=\mathbb{Z}[\overline{y}_{{4}%
}]/\left\langle h_{{8}},h_{{12}}\right\rangle $.
\end{enumerate}

\textbf{Step 3. }The proof is completed by $d_{{31}}=\pm \overline{s}_{{4,2}%
}d_{{23}}$ (Lemma 9) and $d_{{23}}^{2}\in H^{46}=0$ (see in the first column
of the table).$\square $

\bigskip

\noindent \textbf{Remark 2.} In the ring $\mathbb{Z}[\overline{y}_{{4}}]$
one has $\left\langle h_{{8}},h_{{12}}\right\rangle =\left\langle h_{{8}},26%
\overline{y}_{{4}}^{3}\right\rangle $.$\square $

\bigskip

\noindent \textbf{Theorem 10.}\textsl{\ Let }$y_{{3}}$\textsl{, }$y_{{4}}$%
\textsl{,} $y_{{6}}$\textsl{\ be the Schubert classes on }$E_{{6}}/A_{{6}%
}\cdot S^{1}$ \textsl{with Weyl coordinates }$\sigma \lbrack 5,4,2]$\textsl{%
, }$\sigma \lbrack {6,5,4,2}]$\textsl{, }$\sigma \lbrack 1,3,6,5,4,2]$%
\textsl{\ respectively, and let }$d_{{23}},d_{{29}}\in H^{\text{odd}}(E_{{6}%
}/A_{{6}})$\textsl{\ be with }

$\qquad \beta (d_{{23}})=2s_{{11,1}}-s_{{11,2}}$, $\quad \beta (d_{{29}})=s_{%
{14,1}}+s_{{14,2}}+s_{{14,4}}-\ s_{{14,5}}$.

\noindent \textsl{Then}

\begin{center}
$H^{\ast }(E_{{6}}/A_{{6}})=\{\mathbb{Z}[\overline{y}_{{3}},\overline{y}_{{4}%
},\overline{y}_{{6}}]/\left\langle h_{{6}},h_{{8}},h_{{9}},h_{{12}%
}\right\rangle \widehat{\otimes }\Gamma (1,d_{{23}},d_{{29}})\}/\left\langle
2d_{{29}}=\overline{y}_{{3}}d_{{23}}\right\rangle $,
\end{center}

\noindent \textsl{where} $h_{{6}}=2\overline{y}_{{6}}+\overline{y}_{{3}}^{2}$%
\textsl{, }$h_{{8}}=3\overline{y}_{{4}}^{2}$\textsl{, }$h_{{9}}=2\overline{y}%
_{{3}}\overline{y}_{{6}}$\textsl{, }$h_{{12}}=\overline{y}_{{6}}^{2}-%
\overline{y}_{{4}}^{3}$\textsl{.}

\textbf{Proof. Step 1. }From the matrices\textbf{\ }$A_{{k}}$ presented in
[DZ$_{{2}}$, 3.2], one deduces the results in the first two columns of the
table below.

\begin{center}
\begin{tabular}{|l|l|l|}
\hline
nontrivial $H^{k}$ & basis elements & relations \\ \hline
$H^{6}\cong \mathbb{Z}$ & $\bar{s}_{{3,2}}$ &  \\ \hline
$H^{8}\cong \mathbb{Z}$ & $\bar{s}_{{4,3}}$ &  \\ \hline
$H^{12}\cong \mathbb{Z}$ & $\bar{s}_{{6,1}}$ & $-2\bar{s}_{{6,1}}=\bar{s}_{{%
3,2}}^{2}$ \\ \hline
$H^{14}\cong \mathbb{Z}$ & $\bar{s}_{{7,1}}$ & $\overline{s}_{3,2}\,%
\overline{s}_{4,3}$ \\ \hline
$H^{16}\cong \mathbb{Z}_{{3}}$ & $\bar{s}_{{8,1}}$ & $\overline{{s}}{_{4,3}}%
^{2}$ \\ \hline
$H^{18}\cong \mathbb{Z}_{{2}}$ & $\bar{s}_{{9,1}}$ & $\overline{s}_{3,2}\,%
\overline{s}_{6,1}$ \\ \hline
$H^{20}\cong \mathbb{Z}$ & $\bar{s}_{{10,1}}$ & $-\overline{s}_{4,3}\,%
\overline{s}_{6,1}$ \\ \hline
$H^{22}\cong \mathbb{Z}_{{3}}$ & $\bar{s}_{{11,1}}$ & $\overline{s}_{4,3}^{2}%
\overline{s}_{3,2}$ \\ \hline
$H^{26}\cong \mathbb{Z}_{{2}}$ & $\bar{s}_{{13,2}}$ & $\overline{s}_{3,2}\,%
\overline{s}_{4,3}\overline{s}_{6,1}\,$ \\ \hline
$H^{28}\cong \mathbb{Z}_{{3}}$ & $\bar{s}_{{14,1}}$ & $-\,\overline{s}%
_{4,3}^{2}\overline{s}_{6,1}$ \\ \hline
$H^{23}\cong \mathbb{Z}$ & $d_{{23}}=\beta
^{-1}(s_{11,1}-s_{11,2}-s_{11,3}+s_{11,4}$ &  \\ 
& \ \ \ \ \ \ \ \ $-s_{{11,5}}+s_{{11,6}})$ &  \\ \hline
$H^{29}\cong \mathbb{Z}$ & $d_{{29}}=\beta
^{-1}(-s_{14,1}+s_{14,2}+s_{14,4}-\ s_{14,5})$ & $2d_{{29}}=\pm \bar{s}_{{3,2%
}}d_{{23}}$ \\ \hline
$H^{31}\cong \mathbb{Z}$ & $d_{{31}}=\beta
^{-1}(s_{15,1}-2\,s_{15,2}+s_{15,3}-\ s_{15,4})$ & $\pm \bar{s}_{{4,3}}d_{{23%
}}$ \\ \hline
$H^{35}\cong \mathbb{Z}$ & $d_{{35}}=\beta
^{-1}(-s_{17,1}+s_{17,2}+s_{17,3}) $ & $\pm \bar{s}_{{6,1}}d_{{23}}$ \\ 
\hline
$H^{37}\cong \mathbb{Z}$ & $d_{{37}}=\beta ^{-1}(-s_{18,1}+s_{18,2})$ & $\pm 
\bar{s}_{{4,3}}d_{{29}}$ \\ \hline
$H^{43}\cong \mathbb{Z}$ & $d_{{43}}=\beta ^{-1}(s_{{22,1}})$ & $\pm \bar{s}%
_{{4,3}}\bar{s}_{{6,1}}d_{{23}}$ \\ \hline
\end{tabular}
\end{center}

\textbf{Step 2.} From the results in the second column one finds that a
basis of $H^{even}(E_{{6}}/A_{{6}})$ is given as $p^{\ast }\Lambda $, where $%
\Lambda =\{s_{{3,2}},s_{{4,3}},s_{{6,1}},s_{{7,1}},s_{{8,1}},$ $s_{{9,1}},s_{%
{10,1}},$ $s_{{11,1}},s_{{13,1}},s_{{14,1}}\}$. With the multiplicative rule
(4.12) for these basis elements being determined and printed in [DZ$_{{2}}$,
3.3], the items in the last column corresponding to $H^{even}$ are
established. These imply that, if we put $y_{{3}}=s_{{3,2}},y_{{4}}=s_{{{4,}3%
}}{,}y_{{6}}=s_{{6,1}}$, then

\begin{quote}
a) $y_{{3}}$, $y_{{4}}${, }$y_{{6}}$ are the Schubert classes whose Weyl
coordinates are given as those in the theorem by [DZ$_{{2}}$, 3.1];

b) $H^{\text{even}}(E_{{6}}/A_{{6}})$ is generated by $\overline{y}_{{3}}$, $%
\overline{y}_{{4}}$, $\overline{y}_{{6}}$ subject to the relations $h_{{6}%
},h_{{8}},h_{{9}}$.
\end{quote}

\noindent Combining these with the obvious relations $\overline{y}_{{6}}^{2}=%
\overline{y}_{{4}}{^{3}=0}\in H^{24}=0$, together with the fact that $%
\left\langle h_{{6}},h_{{8}},h_{{9}},\overline{y}_{{6}}^{2},\overline{y}_{{4}%
}{^{3}}\right\rangle {=}\left\langle h_{{6}},h_{{8}},h_{{9}},h_{{12}%
}\right\rangle $ in $\mathbb{Z}[\overline{y}_{{3}},\overline{y}_{{4}}{,}%
\overline{y}_{{6}}]$, yields

\begin{enumerate}
\item[(5.4)] $H^{\text{even}}(E_{{6}}/A_{{6}})=\mathbb{Z}[\overline{y}_{{3}},%
\overline{y}_{{4}},\overline{y}_{{6}}]/\left\langle h_{{6}},h_{{8}},h_{{9}%
},h_{{12}}\right\rangle $.
\end{enumerate}

\textbf{Step 3. }By results in the second column $H^{2k+1}\cong \mathbb{Z}$
is spanned by $d_{{2k+1}}$, $k=11,14,15,17,18,21$. From the results in the
first column we deduce also

\begin{quote}
$d_{{2k+1}}d_{{2k^{\prime }+1}}\in H^{2(k+k^{\prime }+1)}=0$ for $%
k,k^{\prime }\in \{11,14,15,17,18,21\}$.
\end{quote}

\noindent Further, we may assume, by the degree reasons, that

\begin{quote}
$\overline{s}_{{3,2}}d_{{23}}=a_{{1}}d_{{29}}$; $\overline{s}_{{4,3}}d_{{23}%
}=a_{{2}}d_{{31}}$; $\overline{s}_{{6,1}}d_{{23}}=a_{{3}}d_{{35}}$; $%
\overline{s}_{{4,3}}d_{{29}}=a_{{4}}d_{{37}}$.
\end{quote}

Lemma 9 suffices to determine the $a_{{i}}\in \mathbb{Z}$ up to sign. For
instance, applying it to the pairings $H^{20}\times H^{23}\rightarrow H^{43}$%
, $H^{14}\times H^{29}\rightarrow H^{43}$ yields respectively that $d_{{43}%
}=\pm \overline{s}_{{4,3}}\overline{s}_{{6,1}}d_{{23}}$, $d_{{43}}=\pm 
\overline{s}_{{3,2}}\overline{s}_{{4,3}}d_{{29}}$. These imply that

\begin{quote}
$\overline{s}_{{4,3}}\overline{s}_{{6,1}}d_{{23}}=\pm \overline{s}_{{3,2}}%
\overline{s}_{{4,3}}d_{{29}}$

$\qquad \ \ \ \ \ \ \ {}=\pm a_{{1}}^{-1}\overline{s}_{{3,2}}^{2}\overline{s}%
_{{4,3}}d_{{23}}$ (by the assumption $\overline{s}_{{3,2}}d_{{23}}=a_{{1}}d_{%
{29}}$)

$\qquad \ \ \ \ \ \ \ {}=\pm 2a_{{1}}^{-1}\overline{s}_{{4,3}}\overline{s}_{{%
6,1}}d_{{23}}$ (by $h_{{6}}$).
\end{quote}

\noindent Coefficients comparison gives $a_{{1}}=\pm 2$.

The same method is applicable to show $a_{{i}}=\pm 1$, $i=2,3,4$. These
verify the items in the third column of the table corresponding to $H^{\text{%
odd}}$, and therefore, complete the proof of Theorem 10.$\square $

\bigskip\ 

\noindent \textbf{Theorem 11.}\ \textsl{Let }$y_{{4}}$\textsl{\ be the
Schubert class on }$E_{{6}}/D_{{5}}\cdot S^{1}$ \textsl{with Weyl coordinate 
}$\sigma \lbrack 2,4,5,6]$\textsl{, and let }$d_{{17}}\in H^{\text{odd}}(E_{{%
6}}/D_{{5}})$\textsl{\ be with }$\beta (d_{{17}})=s_{{8,1}}-s_{{8,2}}-s_{{8,3%
}}$\textsl{. Then}

\begin{center}
$H^{\ast }(E_{{6}}/D_{{5}})=\mathbb{Z}[\overline{y}_{{4}}]/\left\langle h_{{%
12}}\right\rangle \widehat{\otimes }\Gamma (1,d_{{17}})$\textsl{, }
\end{center}

\noindent \textsl{where }$h_{{12}}=\overline{y}_{{4}}^{3}$\textsl{.}

\noindent \textbf{Proof. }The results in [DZ$_{{2}}$, 4.2, 4.3] are
summarized in the table below.

\begin{center}
\begin{tabular}{|l|l|l|}
\hline
nontrivial $H^{k}$ & basis elements & relations \\ \hline
$H^{8}\cong \mathbb{Z}$ & $\overline{s}_{{4,1}}$ &  \\ \hline
$H^{16}\cong \mathbb{Z}$ & $\overline{s}_{{8,1}}$ & $\overline{s}_{{4,1}%
}^{2} $ \\ \hline
$H^{17}\cong \mathbb{Z}$ & $d_{{17}}=\beta ^{-1}(s_{{8,1}}-s_{{8,2}}-s_{{8,3}%
})$ &  \\ \hline
$H^{25}\cong \mathbb{Z}$ & $d_{{25}}=\beta ^{-1}(s_{{12,1}}-s_{{12,2}})$ & $%
\pm \overline{s}_{{4,1}}d_{{17}}$ \\ \hline
$H^{33}\cong \mathbb{Z}$ & $d_{{33}}=\beta ^{-1}(s_{{16,1}})$ & $\pm 
\overline{s}_{{4,1}}^{2}d_{{17}}$ \\ \hline
\end{tabular}%
.
\end{center}

\noindent These imply that

\begin{enumerate}
\item[(5.5)] $H^{\text{even}}(E_{{6}}/D_{{5}})=\mathbb{Z}[\overline{y}_{{4}%
}]/\left\langle h_{{12}}\right\rangle $,
\end{enumerate}

\noindent where $y_{{4}}=s_{{4,1}}$ is the Schubert class whose Weyl
coordinate is given as that in the theorem by [DZ$_{{2}}$, 4.1].

Since $H^{25}\cong \mathbb{Z}$ is generated $d_{{25}}$, $\overline{s}_{{4,1}%
}d_{{17}}=ad_{{25}}$ for some $a\in \mathbb{Z}$. Applying Lemma 9 to the
pairings $H^{8}\times H^{25}\rightarrow H^{33}$, $H^{16}\times
H^{17}\rightarrow H^{33}$ yield respectively that $d_{{33}}=\pm \overline{s}%
_{{4,1}}d_{{25}}$, $d_{{33}}=\pm \overline{s}_{{4,1}}^{2}d_{{17}}$. These
imply that $a=\pm 1$. The proof is completed by $d_{{17}}^{2}\in H^{34}=0$
(see in the first column of the table).$\square $

\bigskip

\noindent \textbf{Theorem 12.}\ \textsl{Let }$y_{{5}}$\textsl{,} $y_{{9}}$%
\textsl{\ be the Schubert classes on }$E_{{7}}/E_{{6}}\cdot S^{1}$ \textsl{%
with Weyl coordinates }$\sigma \lbrack 2,4,5,6,7]$\textsl{, }$\sigma \lbrack
1,5,4,2,3,4,5,6,7]$\textsl{\ respectively, and let }$d_{{37}}$\textsl{, }$d_{%
{45}}\in H^{\text{odd}}(E_{{7}}/E_{{6}})$\textsl{\ be with }$\beta (d_{{37}%
})=s_{{18,1}}-s_{{18,2}}+s_{{18,3}}$\textsl{, }$\beta (d_{{45}})=s_{{22,1}%
}-s_{{22,2}}$\textsl{. Then}

\begin{center}
$H^{\ast }(E_{{7}}/E_{{6}})=\{\mathbb{Z}[\overline{y}_{{5}},\overline{y}_{{9}%
}]/\left\langle h_{{10}},h_{{14}},h_{{18}}\right\rangle \widehat{\otimes }%
\Gamma (1,d_{{37}},d_{{45}})\}/\left\langle \overline{y}_{{9}}d_{{37}}=%
\overline{y}_{{5}}d_{{45}}\right\rangle $\textsl{,}
\end{center}

\noindent \textsl{where }$h_{{10}}=\overline{y}_{{5}}^{2}$\textsl{;} $h_{{14}%
}=2\overline{y}_{{5}}\overline{y}_{{9}}$\textsl{;} $h_{{18}}=\overline{y}_{{9%
}}^{2}$\textsl{.}

\textbf{Proof. }The results in [DZ$_{{2}}$, 5.2, 5.3] are summarized in the
table below.

\begin{center}
\begin{tabular}{|l|l|l|}
\hline
nontrivial $H^{k}$ & basis elements & relations \\ \hline
$H^{10}\cong \mathbb{Z}$ & $\overline{s}_{{5,1}}$ & $\overline{s}_{{5,1}}$
\\ \hline
$H^{18}\cong \mathbb{Z}$ & $\overline{s}_{{9,1}}$ & $\overline{s}_{{9,1}}$
\\ \hline
$H^{28}\cong \mathbb{Z}_{{2}}$ & $\overline{s}_{{14,1}}$ & $\overline{s}_{{%
5,1}}\overline{s}_{{9,1}}$ \\ \hline
$H^{37}\cong \mathbb{Z}$ & $d_{{37}}=\beta ^{-1}(s_{{18,1}}-s_{{18,2}}+s_{{%
18,3}})$ &  \\ \hline
$H^{45}\cong \mathbb{Z}$ & $d_{{45}}=\beta ^{-1}(s_{{22,1}}-s_{{22,2}})$ & 
\\ \hline
$H^{55}\cong \mathbb{Z}$ & $d_{{55}}=\beta ^{-1}(s_{{27,1}})$ & $\overline{s}%
_{{9,1}}d_{{37}}=\pm \overline{s}_{{5,1}}d_{{45}}$ \\ \hline
\end{tabular}
\end{center}

\noindent These imply that, if we let $y_{{5}}=s_{{5,1}}$, $y_{{9}}=s_{{9,1}%
} $, then $y_{{5}}$, $y_{{9}}$ are the Schubert classes whose Weyl
coordinates are given as those in the theorem by [DZ$_{{2}}$, 5.1], and

\begin{enumerate}
\item[(5.6)] $H^{\text{even}}(E_{{7}}/E_{{6}})=\mathbb{Z}[\overline{y}_{{5}},%
\overline{y}_{{9}}]/\left\langle h_{{10}},h_{{14}},h_{{18}}\right\rangle $.
\end{enumerate}

Applying Lemma 9 to the pairings $H^{10}\times H^{45}\rightarrow H^{55}$, $%
H^{18}\times H^{37}\rightarrow H^{55}$ yields that $\overline{s}_{{9,1}}d_{{%
37}}=\pm \overline{s}_{{5,1}}d_{{45}}$. In addition, from $H^{74}=H^{90}=0$
by first column of the table, one gets $d_{{37}}^{2}=$ $d_{{45}}^{2}=0$.
These complete the proof.$\square $

\bigskip

\noindent \textbf{Theorem 13.}\ \textsl{Let }$y_{{4}},y_{{6}},y_{{9}}$ 
\textsl{be the Schubert classes on }$E_{{7}}/D_{{6}}\cdot S^{1}$ \textsl{%
with Weyl coordinates }$\sigma \lbrack 2,4,3,1]$\textsl{, }$\sigma \lbrack
2,6,5,4,3,1]$\textsl{, }$\sigma \lbrack 3,4,2,7,6,5,4,3,1]$\textsl{\
respectively, and let }$d_{{35}},d_{{51}}\in H^{\text{odd}}(E_{{7}}/D_{{6}})$%
\textsl{\ be with}

$\qquad \beta (d_{{35}})=s_{{17,1}}-s_{{17,2}}-s_{{17,3}}+s_{{17,4}}-s_{{17,5%
}}+s_{{17,6}}-s_{{17,7}}$\textsl{;}

$\qquad \beta (d_{{51}})=s_{{25,1}}-s_{{25,2}}-s_{{25,4}}$\textsl{.}

\noindent \textsl{Then}

\begin{center}
$H^{\ast }(E_{{7}}/D_{{6}})=\{\mathbb{Z}[\overline{y}_{{4}},\overline{y}_{{6}%
},\overline{y}_{{9}}]/\left\langle h_{{9}},h_{{12}},h_{{14}},h_{{18}%
}\right\rangle \widehat{\otimes }\Gamma (1,d_{{35}},d_{{51}})\}/\left\langle
3d_{{51}}=\overline{y}_{{4}}^{2}d_{{35}}\right\rangle $\textsl{,}
\end{center}

\noindent \textsl{where} $h_{{9}}=2\overline{y}_{{9}}$\textsl{, }$h_{{12}}=3%
\overline{y}_{{6}}^{2}-\overline{y}_{{4}}^{3}$\textsl{, }$h_{{14}}=3%
\overline{y}_{{4}}^{2}\overline{y}_{{6}}$\textsl{, }$h_{{18}}=\overline{y}_{{%
9}}^{2}-\overline{y}_{{6}}^{3}$\textsl{.}

\textbf{Proof.} The results in [DZ$_{{2}}$, 6.2, 6.3] are summarized in the
table below.

\begin{center}
\begin{tabular}{|l|l|l|}
\hline
nontrivial $H^{k}$ & basis elements & relations \\ \hline
$H^{8}\cong \mathbb{Z}$ & $\overline{s}_{{4,1}}$ &  \\ \hline
$H^{12}\cong \mathbb{Z}$ & $\overline{s}_{{6,1}}$ &  \\ \hline
$H^{16}\cong \mathbb{Z}$ & $\overline{s}_{{8,1}}$ & $\overline{s}_{{4,1}%
}^{2} $ \\ \hline
$H^{18}\cong \mathbb{Z}_{{2}}$ & $\overline{s}_{{9,2}}$ &  \\ \hline
$H^{20}\cong \mathbb{Z}$ & $\overline{s}_{{10,1}}$ & $\overline{s}_{{4,1}}%
\overline{s}_{{6,1}}$ \\ \hline
$H^{24}\cong \mathbb{Z}$ & $\overline{s}_{{12,2}}$ & $\overline{s}_{{12,2}}=%
\overline{s}_{{6,1}}^{2};3\overline{s}_{{12,2}}=\overline{s}_{{4,1}}^{3}$ \\ 
\hline
$H^{26}\cong \mathbb{Z}_{{2}}$ & $\overline{s}_{{13,1}}$ & $\overline{s}_{{%
4,1}}\overline{s}_{{9,2}}$ \\ \hline
$H^{28}\cong \mathbb{Z}_{{3}}$ & $\overline{s}_{{14,1}}$ & $-\overline{s}_{{%
4,1}}^{2}\overline{s}_{{6,1}}$ \\ \hline
$H^{30}\cong \mathbb{Z}_{{2}}$ & $\overline{s}_{{15,1}}$ & $\overline{s}_{{%
6,1}}\overline{s}_{{9,2}}$ \\ \hline
$H^{32}\cong \mathbb{Z}$ & $\overline{s}_{{16,1}}$ & $\overline{s}_{{4,1}}%
\overline{s}_{{6,1}}^{2}$ \\ \hline
$H^{34}\cong \mathbb{Z}_{{2}}$ & $\overline{s}_{{17,2}}$ & $\overline{s}_{{%
4,1}}^{2}\overline{s}_{{9,2}}$ \\ \hline
$H^{38}\cong \mathbb{Z}_{{2}}$ & $\overline{s}_{{19,2}}$ & $\overline{s}_{{%
4,1}}\overline{s}_{{6,1}}\overline{s}_{{9,2}}$ \\ \hline
$H^{40}\cong \mathbb{Z}_{{3}}$ & $\overline{s}_{{20,1}}$ & $\overline{s}_{{%
4,1}}^{2}\overline{s}_{{6,1}}^{2}$ \\ \hline
$H^{42}\cong \mathbb{Z}_{{2}}$ & $\overline{s}_{{21,3}}$ & $\overline{s}_{{%
4,1}}^{3}\overline{s}_{{9,2}}$ \\ \hline
$H^{50}\cong \mathbb{Z}_{{2}}$ & $\overline{s}_{{25,1}}$ & $\overline{s}_{{%
4,1}}^{4}\overline{s}_{{9,2}}$ \\ \hline
$H^{35}\cong \mathbb{Z}$ & $d_{{35}}=\beta ^{-1}(s_{{17,1}}-s_{{17,2}}-s_{{%
17,3}}$ &  \\ 
& \ \ \ $+s_{{17,4}}-s_{{17,5}}+s_{{17,6}}-s_{{17,7}})$ &  \\ \hline
$H^{43}\cong \mathbb{Z}$ & $\beta ^{-1}(s_{{21,1}}-2\,s_{{21,2}}+s_{{21,3}}$
& $\pm \overline{s}_{{4,1}}d_{{35}}$ \\ 
& \ \ \ $-3\,s_{{21,4}}+2\,s_{{21,5}}-s_{{21,6}})$ &  \\ \hline
$H^{47}\cong \mathbb{Z}$ & $\beta ^{-1}(2\,s_{{23,1}}-s_{{23,2}}+s_{{23,3}%
}-s_{{23,4}}$ & $\pm \overline{s}_{{6,1}}d_{{35}}$ \\ 
& \qquad $+s_{{23,5}})$ &  \\ \hline
$H^{51}\cong \mathbb{Z}$ & $d_{{51}}=\beta ^{-1}(s_{{25,1}}-s_{{25,2}}-s_{{%
25,4}})$ & $3d_{{51}}=\pm \overline{s}_{{4,1}}^{2}d_{{35}}$ \\ \hline
$H^{55}\cong \mathbb{Z}$ & $\beta ^{-1}(s_{{27,1}}+s_{{27,2}}-s_{{27,3}})$ & 
$\pm \overline{s}_{{4,1}}\overline{s}_{{6,1}}d_{{35}}$ \\ \hline
$H^{59}\cong \mathbb{Z}$ & $\beta ^{-1}(s_{{29,1}}-s_{{29,2}})$ & $\pm 
\overline{s}_{{6,1}}^{2}d_{{35}},\pm \overline{s}_{{4,1}}d_{{51}}$ \\ \hline
$H^{67}\cong \mathbb{Z}$ & $\beta ^{-1}(s_{{33,1}})$ & $\overline{s}_{{4,1}}%
\overline{s}_{{6,1}}^{2}\!d_{{35}\!}=\!\pm \overline{s}_{{4,1}}^{2}\!d_{{51}%
} $ \\ \hline
\end{tabular}
\end{center}

\noindent These imply that, if we put $y_{{4}}=s_{{4,1}}$, $y_{{6}}=s_{{6,1}%
} ${, }$y_{{9}}=s_{{9,2}}$, then

\begin{quote}
a) $y_{{4}},y_{{6}}{,}y_{{9}}$ are Schubert classes whose Weyl coordinates
are given as those in the theorem by [DZ$_{{2}}$, 6.1];

b) $H^{\text{even}}(E_{{7}}/D_{{6}})$ is generated by $\overline{y}_{{4}},%
\overline{y}_{{6}}{,}\overline{y}_{{9}}$ subject to the relations $h_{{9}%
},h_{{12}},h_{{14}}$.
\end{quote}

\noindent Combining these with the obvious relations $\overline{y}_{{9}}^{2}=%
\overline{y}_{{6}}^{3}\in $ $H^{36}=0$ (see in the first column of the
table), together with $\left\langle h_{{9}},h_{{12}},h_{{14}},\overline{y}_{{%
9}}^{2},\overline{y}_{{6}}^{3}\right\rangle {=}\left\langle h_{{9}},h_{{12}%
},h_{{14}},h_{{18}}\right\rangle $ in $\mathbb{Z}[\overline{y}_{{4}},%
\overline{y}_{{6}},\overline{y}_{{9}}]$, one obtains

\begin{enumerate}
\item[(5.7)] $H^{\text{even}}(E_{{7}}/D_{{6}})=\mathbb{Z}[\overline{y}_{{4}},%
\overline{y}_{{6}},\overline{y}_{{9}}]/\left\langle h_{{9}},h_{{12}},h_{{14}%
},h_{{18}}\right\rangle $.
\end{enumerate}

Finally, the same argument as those in Step 3 in the proof of Theorem 9
verifies the items in the last column of the table corresponding to $H^{%
\text{odd}}$.$\square $

\bigskip

\noindent \textbf{Remark 3.} Let $h_{{9}},h_{{12}},h_{{14}},h_{{18}}$ be the
polynomials in Theorem 12. It can be shown that $\left\langle h_{{9}},h_{{12}%
},h_{{14}},h_{{18}}\right\rangle =\left\langle h_{{9}},h_{{12}},h_{{14}},5%
\overline{y}_{{9}}^{2}+29\overline{y}_{{6}}^{3}\right\rangle $ in $\mathbb{Z}%
[\overline{y}_{{4}},\overline{y}_{{6}},\overline{y}_{{9}}]$.$\square $

\bigskip

\noindent \textbf{Theorem 14.}\ \textsl{Let }$y_{{6}}$\textsl{,}$y_{{10}},y_{%
{15}}$\textsl{\ be the Schubert classes on }$E_{{8}}/E_{{7}}\cdot S^{1}$ 
\textsl{with Weyl coordinates\ }$\sigma \lbrack 3,4,5,6,7,8]$, $\sigma
\lbrack 1,5,4,2,3,4,5,6,7,8]$, $\sigma \lbrack 5,4,3,1,7,6,5,4,2,3,4,$

\noindent $5,6,7,8]$ \textsl{respectively, and let }$d_{{59}}\in H^{\text{odd%
}}(E_{{8}}/E_{{7}})$\textsl{\ be with}

\textsl{\ }$\beta (d_{{59}})=s_{{29,1}}-s_{{29,2}}-s_{{29,3}}+\ s_{{29,4}%
}-s_{{29,5}}+s_{{29,6}}-s_{{29,7}}+s_{{29,8}}$\textsl{. }

\noindent \textsl{Then}

\begin{center}
$H^{\ast }(E_{{8}}/E_{{7}})=\mathbb{Z}[\overline{y}_{{6}},\overline{y}_{{10}%
},\overline{y}_{{15}}]/\left\langle h_{{15}},h_{{20}},h_{{24}},h_{{30}%
}\right\rangle \widehat{\otimes }\Gamma (1,d_{{59}})$\textsl{,}
\end{center}

\noindent \textsl{where }$h_{{15}}=2\overline{y}_{{15}}$\textsl{,} $h_{{20}%
}=3\overline{y}_{{10}}^{2}$\textsl{,} $h_{{24}}=5\overline{y}_{{6}}^{4}$%
\textsl{, }$h_{{30}}=\overline{y}_{{6}}^{5}+\overline{y}_{{10}}^{3}+%
\overline{y}_{{15}}^{2}=0$.

\textbf{Proof. }The results in [DZ$_{{2}}$, 7.2, 7.3] are summarized in the
table below.

\begin{center}
\begin{tabular}{|l|l|l|}
\hline
nontrivial $H^{k}$ & basis elements & relations \\ \hline
$H^{12}\cong \mathbb{Z}$ & $\bar{s}_{{6,2}}$ &  \\ \hline
$H^{20}\cong \mathbb{Z}$ & $\bar{s}_{{10,1}}$ &  \\ \hline
$H^{24}\cong \mathbb{Z}$ & $\bar{s}_{{12,1}}$ & $\pm \bar{s}_{{6,2}}^{2}$ \\ 
\hline
$H^{30}\cong \mathbb{Z}_{{2}}$ & $\bar{s}_{{15,4}}$ &  \\ \hline
$H^{32}\cong \mathbb{Z}$ & $\bar{s}_{{16,1}}$ & $\pm \bar{s}_{{6,2}}\bar{s}_{%
{10,1}}$ \\ \hline
$H^{36}\cong \mathbb{Z}$ & $\bar{s}_{{18,2}}$ & $\pm \bar{s}_{{6,2}}^{3}$ \\ 
\hline
$H^{40}\cong \mathbb{Z}_{{3}}$ & $\bar{s}_{{20,1}}$ & $\pm \bar{s}_{{10,1}%
}^{2}$ \\ \hline
$H^{42}\cong \mathbb{Z}_{{2}}$ & $\bar{s}_{{21,3}}$ & $\pm \bar{s}_{{6,2}}%
\bar{s}_{{15,4}}$ \\ \hline
$H^{44}\cong \mathbb{Z}$ & $\bar{s}_{{22,1}}$ & $\pm \bar{s}_{{6,2}}^{2}\bar{%
s}^{{10,1}}$ \\ \hline
$H^{48}\cong \mathbb{Z}_{{5}}$ & $\bar{s}_{{24,1}}$ & $\pm \bar{s}_{{6,2}%
}^{4}$ \\ \hline
$H^{50}\cong \mathbb{Z}_{{2}}$ & $\bar{s}_{{25,1}}$ & $\pm \bar{s}_{{10,1}}%
\bar{s}_{{15,4}}$ \\ \hline
$H^{52}\cong \mathbb{Z}_{{3}}$ & $\bar{s}_{{26,1}}$ & $\pm \bar{s}_{{6,2}}%
\bar{s}_{{10,1}}^{2}$ \\ \hline
$H^{54}\cong \mathbb{Z}_{{2}}$ & $\bar{s}_{{27,1}}$ & $\pm \bar{s}_{{6,2}%
}^{2}\bar{s}_{{15,4}}$ \\ \hline
$H^{56}\cong \mathbb{Z}$ & $\bar{s}_{{28,1}}$ & $\pm \bar{s}_{{6,2}}^{3}\bar{%
s}_{{10,1}}$ \\ \hline
$H^{62}\cong \mathbb{Z}_{{2}}$ & $\bar{s}_{{31,2}}$ & $\pm \bar{s}_{{6,2}}%
\bar{s}_{{10,1}}\bar{s}_{{15,4}}$ \\ \hline
$H^{64}\cong \mathbb{Z}_{{3}}$ & $\bar{s}_{{32,1}}$ & $\pm \bar{s}_{{6,2}%
}^{2}\bar{s}_{{10,1}}^{2}$ \\ \hline
$H^{66}\cong \mathbb{Z}_{{2}}$ & $\bar{s}_{{33,3}}$ & $\pm \bar{s}_{{6,2}%
}^{3}\bar{s}_{{15,4}}$ \\ \hline
$H^{68}\cong \mathbb{Z}_{{5}}$ & $\bar{s}_{{34,1}}$ & $\pm \bar{s}_{{6,2}%
}^{4}\bar{s}_{{10,1}}$ \\ \hline
$H^{74}\cong \mathbb{Z}_{{2}}$ & $\bar{s}_{{37,2}}$ & $\pm \bar{s}_{{6,2}%
}^{2}\bar{s}_{{10,1}}\bar{s}_{{15,4}}$ \\ \hline
$H^{76}\cong \mathbb{Z}_{{3}}$ & $\bar{s}_{{38,1}}$ & $\pm \bar{s}_{{6,2}%
}^{3}\bar{s}_{{10,1}}^{2}$ \\ \hline
$H^{86}\cong \mathbb{Z}_{{2}}$ & $\bar{s}_{{43,1}}$ & $\pm \bar{s}_{{6,2}%
}^{3}\bar{s}_{{10,1}}^{2}\bar{s}_{{15,4}}$ \\ \hline
$H^{59}\cong \mathbb{Z}$ & $d_{{59}}=\beta ^{-1}(s_{{29,1}}-s_{{29,2}}-s_{{%
29,3}}+\ s_{{29,4}}$ &  \\ 
& $\ \ \ \ \ -s_{{29,5}}+s_{{29,6}}-s_{{29,7}}+s_{{29,8}})$ &  \\ \hline
$H^{71}\cong \mathbb{Z}$ & $\beta ^{-1}(2s_{{35,1}}-3s_{{35,2}}-\ s_{{35,3}%
}+s_{{35,4}}$ & $\pm \bar{s}_{{6,2}}d_{{59}}$ \\ 
& $\ \ \ \ \ +s_{{35,5}}-s_{{35,6}}+s_{{35,7}})$ &  \\ \hline
$H^{79}\cong \mathbb{Z}$ & $\beta ^{-1}(2s_{{39,1}}-s_{{39,2}}-s_{{39,3}}\
-s_{{39,4}}$ & $\pm \bar{s}_{{10,1}}d_{{59}}$ \\ 
& $\ \ \ \ \ +s_{{39,5}}-2s_{{39,6}})$ &  \\ \hline
$H^{83}\cong \mathbb{Z}$ & $\beta ^{-1}(2\,s_{{41,1}}-s_{{41,2}}+s_{{41,3}%
}-s_{{41,4}}+s_{{41,5}})$ & $\pm \bar{s}_{{6,2}}d_{{59}}$ \\ \hline
$H^{91}\cong \mathbb{Z}$ & $\beta ^{-1}(s_{{45,1}}-s_{{45,2}}-s_{{45,3}}+s_{{%
45,4}})$ & $\pm \bar{s}_{{6,2}}\bar{s}_{{10,1}}d_{{59}}$ \\ \hline
$H^{95}\cong \mathbb{Z}$ & $\beta ^{-1}(s_{{47,1}}-s_{{47,2}}+s_{{47,3}})$ & 
$\pm \bar{s}_{{6,2}}^{3}d_{{59}}$ \\ \hline
$H^{103}\cong \mathbb{Z}$ & $\beta ^{-1}(-s_{{51,1}}+s_{{51,2}})$ & $\pm 
\bar{s}_{{6,2}}^{2}\bar{s}_{{10,1}}d_{{59}}$ \\ \hline
\end{tabular}
\end{center}

\noindent These imply that, if we put $y_{{6}}=s_{{6,2}}$, $y_{{10}}=s_{{10,1%
}},y_{{15}}=s_{{15,4}}$, then

\begin{quote}
a) $y_{{6}},y_{{10}},y_{{15}}$ are the Schubert classes on $E_{{8}}/E_{{7}%
}\cdot S^{1}$ whose Weyl coordinates are given as those in the theorem by [DZ%
$_{{2}}$, 7.1];

b) the ring $H^{\text{even}}(E_{{8}}/E_{{7}})$ is generated by $\overline{y}%
_{{6}},\overline{y}_{{10}},\overline{y}_{{15}}$ subject to the relation $h_{{%
15}}$, $h_{{20}},h_{{24}}$.
\end{quote}

\noindent Combining this with $\overline{y}_{{6}}^{5},\overline{y}_{{10}%
}^{3},\overline{y}_{{15}}^{2}\in H^{60}=0$ (see in the first column), and
noting that $\left\langle h_{{15}},h_{{20}},h_{{24}},h_{{30}}\right\rangle
=\left\langle h_{{15}},h_{{20}},h_{{24}},\overline{y}_{{6}}^{5},\overline{y}%
_{{10}}^{3},\overline{y}_{{15}}^{2}\right\rangle $ in $\mathbb{Z}[\overline{y%
}_{{6}},\overline{y}_{{10}},\overline{y}_{{15}}]$, we get

\begin{enumerate}
\item[(5.8)] $H^{\text{even}}(E_{{8}}/E_{{7}})=\mathbb{Z}[\overline{y}_{{6}},%
\overline{y}_{{10}},\overline{y}_{{15}}]/\left\langle h_{{15}},h_{{20}},h_{{%
24}},h_{{30}}\right\rangle $.
\end{enumerate}

According to results in the second column of the table corresponding to $H^{%
\text{odd}}$, the $H^{2k+1}\cong \mathbb{Z}$ is generated by $d_{{2k+1}}$
with $k\in \{29,35,39,41,45,47,51\}$. Moreover, Lemma 9 is applicable to
show that (see Step 3 in the proof of Theorem 10):

\begin{quote}
$d_{{71}}=\pm \overline{s}_{{6,2}}d_{{59}}$; $d_{{83}}=\pm \overline{s}_{{6,2%
}}d_{{71}}$; $d_{{95}}=\pm \overline{s}_{{6,2}}d_{{83}}$;

$d_{{79}}=\pm \overline{s}_{{10,1}}d_{{59}}$; $d_{{91}}=\pm \overline{s}_{{%
10,1}}d_{{71}}$; $d_{{103}}=\pm \overline{s}_{{10,1}}d_{{83}}$.
\end{quote}

\noindent These complete the proof. $\square $

\section{Proofs of Theorems 1-7}

Let $G/H$ be one of the Grassmannians concerned in Theorems 1--7. Since the
Chow ring $A^{\ast }(G/H)$ is isomorphic to the cohomology $H^{\ast }(G/H)$
via the cycle map (Remark 1), Lemmas 7 and 8 are directly applicable to
compute $A^{\ast }(G/H)$.

Firstly, comparing Lemma 7 with Theorems 8--14 we find that (for each $G/H$)
a minimal set of generators for $A^{\ast }(G/H)$ is given as those asserted
in Theorems 1--7. Next, applying Lemma 8 to the presentations of $H^{\ast
}(G/H_{s})$ in Theorems 8--14 one finds the initial constraints that the
corresponding relations on $A^{\ast }(G/H)$ satisfy. It remains for us to

\begin{quote}
i) specify the relations that are subject to the constraints; and

ii) dispel those relations that belong to the ideal generated by the lower
degree ones.
\end{quote}

\noindent Lemmas 3 and 4 are functional in implementing the tasks i) and ii)
respectively.

As in Lemma 3 we use $\kappa _{i}$ to denote the $i^{th}$ row in a
null-space $N(\pi _{{m}})$.

\bigskip

\noindent \textbf{Proof of Theorem 1. }Combining Theorem 8 with Lemmas 7 and
8, we get the partial description of $A^{\ast }(F_{{4}}/C_{{3}}\cdot S^{1})$:

\begin{enumerate}
\item[(6.1)] $A^{\ast }(F_{{4}}/C_{{3}}\cdot S^{1})=\mathbb{Z}[y_{{1}},y_{{3}%
},y_{{4}},y_{{6}}]/\left\langle r_{{3}},r_{{6}},r_{{8}},r_{{12}},y_{{1}}g_{{%
11}}\right\rangle $,
\end{enumerate}

\noindent where $y_{{1}},y_{{3}},y_{{4}},y_{{6}}$ are the same Schubert
classes as those in Theorem 1, and where if we let $\pi _{{m}}:\mathbb{Z}[y_{%
{1}},y_{{3}},y_{{4}},y_{{6}}]^{(2m)}\rightarrow A^{\ast }(F_{{4}}/C_{{3}%
}\cdot S^{1})$ be induced from $\{y_{{1}},y_{{3}},y_{{4}},y_{{6}}\}\subset
A^{\ast }(F_{{4}}/C_{{3}}\cdot S^{1})$, then (Lemma 8)

\begin{quote}
1) for $m=3,6,8,12$, $r_{{m}}\in \ker \pi _{{m}}$ in (6.1) should satisfy
\end{quote}

$\qquad \qquad r_{{3}}\mid _{{y}_{{1}}{=0}}$ $=2y_{{3}}$; $r_{{6}}\mid _{{y}%
_{{1}}{=0}}$ $=2y_{{6}}+y_{{3}}^{2}$;

$\qquad \qquad r_{{8}}\mid _{{y}_{{1}}{=0}}$ $=3y_{{4}}^{2}$; $r_{{12}}\mid
_{{y}_{{1}}{=0}}$ $=y_{{6}}^{2}-y_{{4}}^{3}$;

\begin{quote}
2) $\pi (g_{{11}})=2\,s_{{11,1}}-s_{{11,2}}$.
\end{quote}

With respect to the ordered basis $B(2m)$ of $\mathbb{Z}[y_{{1}},y_{{3}},y_{{%
4}},y_{{6}}]^{(2m)}$ for $m=3$, $6$, $8$, $12$, the structure matrices $%
M(\pi _{{m}})$ have been computed from the \textsl{L--R coefficients}, and
their corresponding Nullspaces $N(\pi _{{m}})$ are presented in [DZ$_{2}$,
1.5]. If we take, in terms of Lemma 3, that

\begin{quote}
$r_{{3}}=2y_{{3}}-y_{{1}}^{3}$ ($=-\kappa _{{1}}$ in $N(\pi _{{3}})$);

$r_{{6}}=2y_{{6}}+y_{{3}}^{2}-3y_{{1}}^{2}y_{{4}}$ ($=-\kappa _{{3}}$ in $%
N(\pi _{{6}})$);

$r_{{8}}=3y_{{4}}^{2}-y_{{1}}^{2}y_{{6}}$ ($=-\kappa _{{5}}$ in $N(\pi _{{8}%
})$);

$r_{{12}}=y_{{6}}^{2}-y_{{4}}^{3}$ ($=\kappa _{{15}}$ in $N(\pi _{{12}})$),
\end{quote}

\noindent then condition 1) is met by the set $\{r_{{3}},r_{{6}},r_{{8}},r_{{%
12}}\}$ of polynomials.

The proof will be completed once we show

\begin{enumerate}
\item[(6.2)] $y_{{1}}g_{{11}}\in \left\langle r_{{3}},r_{{6}},r_{{8}},r_{{12}%
}\right\rangle $.
\end{enumerate}

\noindent For this purpose we examine, in view of (6.1), the quotient map (%
\S 3.3)

\begin{quote}
$\varphi :\mathbb{Z}[y_{{1}},y_{{3}},y_{{4}},y_{{6}}]/\left\langle r_{{3}%
},r_{{6}},r_{{8}},r_{{12}}\right\rangle \rightarrow A^{\ast }(F_{{4}}/C_{{3}%
}\cdot S^{1})=\tbigoplus_{m\geq 0}A^{m}$.
\end{quote}

\noindent With $r_{{3}},r_{{6}},r_{{8}},r_{{12}}$ being explicitly
presented, it is straightforward to find that

\begin{quote}
$b(24)=16$; $\delta _{{24}}(r_{{3}},r_{{6}},r_{{8}},r_{{12}})=15$ (see
Example 4).
\end{quote}

\noindent On the other hand, granted with the basis theorem, we read from [DZ%
$_{{2}}$, 1.1] that rank($A^{24}$)$=1$. (6.2) is verified by Lemma 4.$%
\square $

\bigskip

\noindent \textbf{Proof of Theorem 2. }Combining Theorem 9 with Lemmas 7 and
8, we get the partial description of $A^{\ast }(F_{{4}}/B_{{3}}\cdot S^{1})$:

\begin{enumerate}
\item[(6.3)] $A^{\ast }(F_{{4}}/B_{{3}}\cdot S^{1})=\mathbb{Z}[y_{{1}},y_{{4}%
}]/\left\langle r_{{8}},r_{{12}},y_{{1}}g_{{11}}\right\rangle $,
\end{enumerate}

\noindent where the generators $y_{{1}},y_{{4}}$ are the same Schubert
classes as those in Theorem 2, and where if we let $\pi _{{m}}:\mathbb{Z}[y_{%
{1}},y_{{4}}]^{(2m)}\rightarrow A^{m}(F_{{4}}/B_{{3}}\cdot S^{1})$ be
induced from $\{y_{{1}},y_{{4}}\}\subset A^{\ast }(F_{{4}}/B_{{3}}\cdot
S^{1})$, then (Lemma 8)

\begin{quote}
1) for $m=8,12$, the $r_{{m}}\in \ker \pi _{{m}}$ in (6.3) should satisfy

$\qquad r_{{8}}\mid _{{y}_{{1}}{=0}}$ $=3y_{{4}}^{2}$; $r_{{12}}\mid _{{y}_{{%
1}}{=0}}$ $=26y_{{4}}^{3}$
\end{quote}

\noindent (see Remark 2 after the proof of Theorem 9);

\begin{quote}
2) $\pi (g_{{11}})=-s_{{11,1}}+s_{{11,2}}$.
\end{quote}

With respect to the ordered basis $B(2m)$ of $\mathbb{Z}[y_{{1}},y_{{4}%
}]^{(2m)}$, $m=8,12$, the structure matrices $M(\pi _{{m}})$ have been
computed by the \textsl{L--R coefficients}, and their corresponding
Nullspaces $N(\pi _{{m}})$ are presented in [DZ$_{2}$, 2.5]. If we take, in
terms of Lemma 3, that

\begin{quote}
$r_{{8}}=3y_{{4}}^{2}-y_{{1}}^{8}$ ($=-\kappa _{{1}}$ in $N(\pi _{{8}})$);

$r_{{12}}=26y_{{4}}^{3}-5y_{{1}}^{12}$ ($=-\kappa _{{1}}$ in $N(\pi _{{12}})$%
),
\end{quote}

\noindent then condition 1) is met by the set $\{r_{{8}},r_{{12}}\}$ of
polynomials.

The proof will be completed once we show

\begin{enumerate}
\item[(6.4)] $y_{{1}}g_{{11}}\in \left\langle r_{{8}},r_{{12}}\right\rangle $%
.
\end{enumerate}

\noindent For this purpose we examine, in view of (6.3), the quotient map (%
\S 3.3)

\begin{center}
$\varphi :\mathbb{Z}[y_{{1}},y_{{4}}]/\left\langle r_{{8}},r_{{12}%
}\right\rangle \rightarrow A^{\ast }(F_{{4}}/B_{{3}}\cdot
S^{1})=\tbigoplus_{m\geq 0}A^{m}$.
\end{center}

\noindent With $r_{{8}},r_{{12}}$ being given explicitly, it is
straightforward to find that

\begin{quote}
$b(24)=4$, $\delta _{{24}}(r_{{8}},r_{{12}})=3$ (see Example 4).
\end{quote}

\noindent On the other hand, granted with the basis theorem, we read from [DZ%
$_{{2}}$, 2.1] that rank($A^{24}$)$=1$. (6.4) is verified by Lemma 4.$%
\square $

\bigskip

\noindent \textbf{Proof of Theorem 3. }Combining Theorem 10 with Lemmas 7
and 8, we get the partial description of $A^{\ast }(E_{{6}}/A_{{6}}\cdot
S^{1})$:

\begin{enumerate}
\item[(6.5)] $A^{\ast }(E_{{6}}/A_{{6}}\cdot S^{1})=\mathbb{Z}[y_{{1}},y_{{3}%
},y_{{4}},y_{{6}}]/\left\langle r_{{6}},r_{{8}},r_{{9}},r_{{12}},y_{{1}}g_{{%
11}},y_{{1}}g_{{14}}\right\rangle $,
\end{enumerate}

\noindent where $y_{{1}},y_{{3}},y_{{4}},y_{{6}}$ are the same Schubert
classes as those in Theorem 3, and where if we let $\pi _{{m}}:\mathbb{Z}[y_{%
{1}},y_{{3}},y_{{4}},y_{{6}}]^{(2m)}\rightarrow A^{\ast }(E_{{6}}/A_{{6}%
}\cdot S^{1})$ be induced from $\{y_{{1}},y_{{3}},y_{{4}},y_{{6}}\}\subset
A^{\ast }(E_{{6}}/A_{{6}}\cdot S^{1})$, then (Lemma 8)

\begin{quote}
1) for $m=6,8,9,12$, $r_{{m}}\in \ker \pi _{{m}}$ in (6.5) should satisfy
\end{quote}

$\qquad \qquad r_{{6}}\mid _{{y}_{{1}}{=0}}$ $=2y_{{6}}+y_{{3}}^{2}$; $r_{{8}%
}\mid _{{y}_{{1}}{=0}}$ $=3y_{{4}}^{2}$;

$\qquad \qquad r_{{9}}\mid _{{y}_{{1}}{=0}}$ $=2y_{{3}}y_{{6}}$; $r_{{12}%
}\mid _{{y}_{{1}}{=0}}$ $=y_{{4}}^{3}-y_{{6}}^{2}$;

\begin{quote}
2) $\pi (g_{{11}})=s_{{11,1}}-s_{{11,2}}-s_{{11,3}}+s_{{11,4}}-\ s_{{11,5}%
}+s_{{11,6}},\quad $

$\quad ~\pi (g_{{14}})=s_{{14,1}}+s_{{14,2}}+s_{{14,4}}-\ s_{{14,5}}$.
\end{quote}

With respect to the ordered basis $B(2m)$ of $\mathbb{Z}[y_{{1}},y_{{3}},y_{{%
4}},y_{{6}}]$, $m=6,8,9,12$, the structure matrices $M(\pi _{{m}})$ have
been computed by the \textsl{L--R coefficients} and their corresponding
Nullspaces $N(\pi _{{m}})$ are presented in [DZ$_{2}$, 3.5]. If we take, in
view of Lemma 3, that

\begin{quote}
$r_{6}$: \ $2{y_{6}}+{y_{3}^{2}}-3{y_{1}}^{2}{y_{4}}+2{y_{1}}^{3}{y_{3}}-{%
y_{1}}^{6}$ ($=\kappa _{1}$ in $N(\pi _{6})$);

$r_{8}$: \ $3{y_{4}}^{2}-6{y_{1}}{y_{3}}{y_{4}}+\ {y_{1}}^{2}{y_{6}}+5{y_{1}}%
^{2}{y_{3}}^{2}-2{y_{1}}^{5}{y_{3}}$ ($=\kappa _{2}$ in $N(\pi _{8})$);

$r_{9}$: \ $2{y_{3}}{y_{6}}-{y_{1}}^{3}{y_{6}}$ ($=\kappa _{4}$ in $N(\pi
_{9})$);

$r_{12}$: \ ${y_{6}}^{2}-{y_{4}}^{3}$ ($=\kappa _{11}$ in $N(\pi _{12})$),
\end{quote}

\noindent then condition 1) is met by the set $\{r_{{6}},r_{{8}},r_{{9}},r_{{%
12}}\}$ of polynomials.

The proof will be completed once we show

\begin{enumerate}
\item[(6.6)] $y_{{1}}g_{{11}},y_{{1}}g_{{14}}\in \left\langle r_{{6}},r_{{8}%
},r_{{9}},r_{{12}}\right\rangle $.
\end{enumerate}

\noindent For this purpose we examine, in view of (6.5), the quotient map
(see\textbf{\ \S }3.3)

\begin{quote}
$\varphi :\mathbb{Z}[y_{{1}},y_{{3}},y_{{4}},y_{{6}}]/\left\langle r_{{6}%
},r_{{8}},r_{{9}},r_{{12}}\right\rangle \rightarrow A^{\ast }(E_{{6}}/A_{{6}%
}\cdot S^{1})=\oplus _{m\geq 0}A^{m}$
\end{quote}

\noindent With $r_{{6}},r_{{8}},r_{{9}},r_{{12}}$ being made explicitly it
is straightforward to find that

\begin{quote}
$b(24)=16$; $\delta _{{24}}(r_{{6}},r_{{8}},r_{{9}},r_{{12}})=11$;

$b(30)=24$; $\delta _{{30}}(r_{{6}},r_{{8}},r_{{9}},r_{{12}})=20$ (Example
4).
\end{quote}

\noindent On the other hand, granted with the basis theorem, we read from [DZ%
$_{{2}}$, 3.1] that rank($A^{24}$)$=5$, rank($A^{30}$)$=4$. (6.6) is
verified by Lemma 4.$\square $

\bigskip

\noindent \textbf{Proof of Theorem 4. }Combining Theorem 11 with Lemmas 7
and 8, we get the partial description for $A^{\ast }(E_{{6}}/D_{{5}}\cdot
S^{1})$:

\begin{enumerate}
\item[(6.7)] $A^{\ast }(E_{{6}}/D_{{5}}\cdot S^{1})=\mathbb{Z}[y_{{1}},y_{{4}%
}]/\left\langle y_{{1}}g_{{8}},r_{{12}}\right\rangle $
\end{enumerate}

\noindent where $y_{{1}},y_{{4}}$ are the same Schubert classes as those in
Theorem 4, and where if we let $\pi _{{m}}:\mathbb{Z}[y_{{1}},y_{{4}%
}]^{(2m)}\rightarrow A^{\ast }(E_{{6}}/D_{{5}}\cdot S^{1})$ be induced from $%
\{y_{{1}},y_{{4}}\}\subset A^{\ast }(E_{{6}}/D_{{5}}\cdot S^{1})$, then

\begin{quote}
1) the $r_{{12}}\in \ker \pi _{{12}}$ in (6.7) should satisfy $r_{{12}}\mid
_{{y}_{{1}}{=0}}=y_{{4}}^{3}$;

2) $\pi (g_{{8}})=s_{{8,1}}-s_{{8,2}}-s_{{8,3}}$.
\end{quote}

Let us find $g_{{8}}\in \mathbb{Z}[y_{{1}},y_{{4}}]$ required to specify the
first relation $y_{{1}}g_{{8}}$. Assume, with respect to the basis $B(16)$
of $\mathbb{Z}[y_{{1}},y_{{4}}]^{(16)}$, that

\begin{enumerate}
\item[(6.8)] $g_{{8}}=a_{{1}}y_{{1}}^{8}+a_{{2}}y_{{1}}^{4}y_{{4}}+a_{{3}}y_{%
{4}}^{2}$.
\end{enumerate}

\noindent According to [DZ$_{{2}}$, 4.1], there are three Schubert classes
in dimension $16$ with Weyl coordinates

$w_{{8,1}}=\sigma \lbrack {1,5,4,2,3,4,5,6}]$; $w_{{8,2}}=\sigma \lbrack {%
3,1,4,2,3,4,5,6}]$;

$w_{{8,3}}=\sigma \lbrack {6,5,4,2,3,4,5,6}]$

\noindent respectively. The constraint 2) implies that $g_{{8}}$ must
satisfy the system

\begin{quote}
$a_{{w}_{{8,1}}}(g_{{8}})=1$; $a_{{w}_{{8,2}}}(g_{{8}})=-1$; $a_{{w}_{{8,3}%
}}(g_{{8}})=1$.
\end{quote}

\noindent Thus, applying the \textsl{L--R Coefficients} (\S 2.3) to (6.8)
yields

\begin{quote}
$\left\{ 
\begin{array}{c}
1=7a_{{1}}+3a_{{2}}+a_{{3}} \\ 
-1=5a_{{1}}+2a_{{2}}+a_{{3}} \\ 
-1=2a_{{1}}+a_{{2}}+a_{{3}}\text{.}%
\end{array}%
\right. $
\end{quote}

\noindent From this we find that $(a_{{1}},a_{{2}},a_{{3}})=(-2,6,-3)$, and
consequently

\begin{quote}
$y_{{1}}g_{{8}}=2y_{{1}}^{9}+3y_{{1}}y_{{4}}^{2}-6y_{{1}}^{5}y_{{4}}$ (see
Theorem 4).
\end{quote}

To find $r_{{12}}$ consider the map $\pi _{{12}}:\mathbb{Z}[y_{{1}},y_{{4}%
}]^{(24)}\rightarrow A^{24}(E_{{6}}/D_{{5}}\cdot S^{1})$. With respect to
the ordered basis $B(24)$ of $\mathbb{Z}[y_{{1}},y_{{4}}]^{(24)}$, the
structure matrix $M(\pi _{{12}})$ has been computed by the \textsl{L--R
coefficients} and the corresponding Nullspaces $N(\pi _{{12}})$ is presented
in [DZ$_{2}$, 4.5]. If we take, in view of Lemma 3, that

\begin{quote}
$r_{{12}}=y_{{4}}^{3}-6y_{{1}}^{4}y_{{4}}^{2}+y_{{1}}^{12}$ ($=\kappa _{{1}}$
in $N(\pi _{{12}})$),
\end{quote}

\noindent then condition 1) is met by the $r_{{12}}$ above. This finishes
the proof.$\square $

\bigskip

\noindent \textbf{Proof of Theorem 5. }Combining Theorem 12 with Lemmas 7
and 8, we get the partial description of $A^{\ast }(E_{{7}}/E_{{6}}\cdot
S^{1})$:

\begin{enumerate}
\item[(6.9)] $A^{\ast }(E_{{7}}/E_{{6}}\cdot S^{1})=\mathbb{Z}[y_{{1}},y_{{5}%
},y_{{9}}]/\left\langle r_{{10}},r_{{14}},r_{{18}},y_{{1}}g_{{18}},y_{{1}}g_{%
{22}}\right\rangle $,
\end{enumerate}

\noindent where $y_{{1}},y_{{5}},y_{{9}}$ are the same Schubert classes as
those in Theorem 5, and where if we let $\pi _{{m}}:\mathbb{Z}[y_{{1}},y_{{5}%
},y_{{9}}]^{(2m)}\rightarrow A^{\ast }(E_{{7}}/E_{{6}}\cdot S^{1})$ be
induced from $\{y_{{1}},y_{{5}},y_{{9}}\}\subset A^{\ast }(E_{{7}}/E_{{6}%
}\cdot S^{1})$, then (Lemma 8)

\begin{quote}
1) for $m=10,14,18$, the $r_{{m}}\in \ker \pi _{{m}}$ in (6.9) should satisfy
\end{quote}

$\ \qquad r_{{10}}\mid _{{y}_{{1}}{=0}}$ $=y_{{5}}^{2}$; $r_{{14}}\mid _{{y}%
_{{1}}{=0}}$ $=2y_{{5}}y_{{9}}$;$\quad r_{{18}}\mid _{{y}_{{1}}{=0}}$ $=y_{{9%
}}^{2}$;

\begin{quote}
2) $\pi (g_{{18}})=s_{{18,1}}-s_{{18,2}}+s_{{18,3}}$,$\quad \pi (g_{{22}})$ $%
=s_{{22,1}}-s_{{22,2}}$.
\end{quote}

With respect to the ordered basis $B(2m)$ of $\mathbb{Z}[y_{{1}},y_{{5}},y_{{%
9}}]^{(2m)}$, $m=10,14,18$, the structure matrices $M(\pi _{{m}})$ have been
computed by the \textsl{L--R coefficients} and their corresponding
Nullspaces $N(\pi _{{m}})$ are presented in [DZ$_{2}$, 5.5]. If we take, in
view of Lemma 3, that

\begin{quote}
$r_{{10}}=y_{{5}}^{2}-2y_{{1}}y_{{9}}$ ($=-\kappa _{{1}}$ in $N(\pi _{{10}})$%
);

$r_{{14}}=2y_{{5}}y_{{9}}-9y_{{1}}^{4}y_{{5}}^{2}+6y_{{1}}^{9}y_{{5}}-y_{{1}%
}^{14}$ ($=-\kappa _{{1}}$ in $N(\pi _{{14}})$);

$r_{{18}}=y_{{9}}^{2}+10y_{{1}}^{3}y_{{5}}^{3}-9y_{{1}}^{8}y_{{5}}^{2}+2y_{{1%
}}^{13}y_{{5}}$ ($=\kappa _{{2}}$ in $N(\pi _{{18}})$),
\end{quote}

\noindent then condition 1) is met by the set $\{r_{{10}},r_{{14}},r_{{18}%
}\} $ of polynomials.

The proof will be completed once we show

\begin{enumerate}
\item[(6.10)] $y_{{1}}g_{{18}}$, $y_{{1}}g_{{22}}\in \left\langle r_{{10}%
},r_{{14}},r_{{18}}\right\rangle $.
\end{enumerate}

\noindent For this purpose we examine, in view of (6.9), the quotient map (%
\S 3.3)

\begin{center}
$\varphi :\mathbb{Z}[y_{{1}},y_{{5}},y_{{9}}]/\left\langle r_{{10}},r_{{14}%
},r_{{18}}\right\rangle \rightarrow A^{\ast }(E_{{7}}/E_{{6}}\cdot
S^{1})=\oplus _{m\geq 0}A^{m}$
\end{center}

\noindent With $r_{{10}},r_{{14}},r_{{18}}$ being obtained explicitly it is
straightforward to find that

\begin{quote}
$b(38)=8$; $\delta _{{38}}(r_{{10}},r_{{14}},r_{{18}})=6$;

$b(46)=10$; $\delta _{{46}}(r_{{10}},r_{{14}},r_{{18}})=9$ (see Example 3).
\end{quote}

\noindent On the other hand, granted with the basis theorem, we read from [DZ%
$_{{2}}$, 5.1] that rank($A^{38}$)$=2$, rank($A^{46}$)$=1$. (6.10) is
verified by Lemma 4.$\square $

\bigskip

\noindent \textbf{Proof of Theorem 6}. Combining Theorem 13 with Lemmas 7
and 8, we get the partial description of $A^{\ast }(E_{{7}}/D_{{6}}\cdot
S^{1})$:

\begin{enumerate}
\item[(6.11)] $A^{\ast }(E_{{7}}/D_{{6}}\cdot S^{1})=\mathbb{Z}[y_{{1}},y_{{4%
}},y_{{6}},y_{{9}}]/\left\langle r_{{9}},r_{{12}},r_{{14}},r_{{18}},y_{{1}%
}g_{{17}},y_{{1}}g_{{25}}\right\rangle $,
\end{enumerate}

\noindent where $y_{{1}},y_{{4}},y_{{6}},y_{{9}}$ are the same Schubert
classes as those in Theorem 6, and where if we let $\pi _{{m}}:\mathbb{Z}[y_{%
{1}},y_{{4}},y_{{6}},y_{{9}}]^{(2m)}\rightarrow A^{\ast }(E_{{7}}/D_{{6}%
}\cdot S^{1})$ be induced from $\{y_{{1}},y_{{4}},y_{{6}},y_{{9}}\}\subset
A^{\ast }(E_{{7}}/D_{{6}}\cdot S^{1})$, then (Lemma 8)

\begin{quote}
1) for $m=9,12,14,18$, the $r_{{m}}\in \ker \pi _{{m}}$ in (6.11) should
satisfy
\end{quote}

$\qquad r_{{9}}\mid _{{y}_{{1}}{=0}}$ $=2y_{{9}}$;$\quad r_{{12}}\mid _{{y}_{%
{1}}{=0}}$ $=3y_{{6}}^{2}-y_{{4}}^{3}$; $r_{{14}}\mid _{{y}_{{1}}{=0}}=3y_{{4%
}}^{2}y_{{6}}$;

$\qquad r_{{18}}\mid _{{y}_{{1}}{=0}}$ $=5y_{{9}}^{2}+29y_{{6}}^{3}$ (Remark
3 after the proof of Theorem 13)

\begin{quote}
2) $\pi (g_{{17}})=s_{{17,1}}-s_{{17,2}}-s_{{17,3}}+s_{{17,4}}-s_{{17,5}}+s_{%
{17,6}}-s_{{17,7}},\quad $

$\ \ \ \pi (g_{{25}})=s_{{25,1}}-s_{{25,2}}-s_{{25,4}}$.
\end{quote}

With respect to the ordered basis $B(2m)$ of $\mathbb{Z}[y_{{1}},y_{{4}},y_{{%
6}},y_{{9}}]^{(2m)}$, $m=9,12,14,$ $18$, the structure matrices $M(\pi _{{m}%
})$ have been computed by the \textsl{L--R coefficients} and their
corresponding Nullspaces $N(\pi _{{m}})$ are presented in [DZ$_{2}$, 6.5].
If we take, in view of Lemma 3, that

\begin{quote}
$r_{{9}}=2y_{{9}}+3y_{{1}}y_{{4}}^{2}+4y_{{1}}^{3}y_{{6}}+2y_{{1}}^{5}y_{{4}%
}-2y_{{1}}^{9}$ ($=-\kappa _{{1}}$ in $N(\pi _{{9}})$);

$r_{{12}}=3y_{{6}}^{2}-y_{{4}}^{3}-3y_{{1}}^{4}y_{{4}}^{2}-2y_{{1}}^{6}y_{{6}%
}+2y_{{1}}^{8}y_{{4}}$($=\kappa _{{1}}$ in $N(\pi _{{12}})$);

$r_{{14}}=3y_{{4}}^{2}y_{{6}}+3y_{{1}}^{2}y_{{6}}^{2}+6y_{{1}}^{2}y_{{4}%
}^{3}+6y_{{1}}^{4}y_{{4}}y_{{6}}+2y_{{1}}^{5}y_{{9}}-y_{{1}}^{14}$

\qquad ($=-\kappa _{{1}}$ in $N(\pi _{{14}})$);

$r_{{18}}=5y_{{9}}^{2}+29y_{{6}}^{3}-24y_{{1}}^{6}y_{{6}}^{2}+45y_{{1}%
}^{2}y_{{4}}y_{{6}}^{2}+2y_{{1}}^{9}y_{{9}}$ ($=\kappa _{{5}}-2\kappa _{{8}}$
in $N(\pi _{{18}})$),
\end{quote}

\noindent then condition 1) is met by the set $\{r_{{9}},r_{{12}},r_{{14}%
},r_{{18}}\}$ of polynomials.

The proof will be completed once we show

\begin{enumerate}
\item[(6.12)] $y_{{1}}g_{{17}},y_{{1}}g_{{25}}\in \left\langle r_{{9}},r_{{12%
}},r_{{14}},r_{{18}}\right\rangle $.
\end{enumerate}

\noindent For this purpose we examine, in view of (6.11), the quotient map (%
\S 3.3)

\begin{center}
$\varphi :\mathbb{Z}[y_{{1}},y_{{4}},y_{{6}},y_{{9}}]/\left\langle r_{{9}%
},r_{{12}},r_{{14}},r_{{18}}\right\rangle \rightarrow A^{\ast }(E_{{7}}/D_{{6%
}}\cdot S^{1})=\oplus _{m\geq 0}A^{m}$
\end{center}

\noindent With $r_{{9}},r_{{12}},r_{{14}},r_{{18}}$ being obtained
explicitly it is straightforward to find that

\begin{quote}
$b(36)=17$; $\delta _{{36}}(r_{{9}},r_{{12}},r_{{14}},r_{{18}})=11$;

$b(52)=32$; $\delta _{{52}}(r_{{9}},r_{{12}},r_{{14}},r_{{18}})=29$ (Example
4).
\end{quote}

\noindent On the other hand, granted with the basis theorem, we read from [DZ%
$_{{2}}$, 6.1] that rank($A^{36}$)$=6$, rank($A^{52}$)$=3$. (6.12) is
verified by Lemma 4.$\square $

\bigskip

\noindent \textbf{Proof of Theorem 7}. Combining Theorem 14 with Lemmas 7
and 8, we get the partial description of $A^{\ast }(E_{{8}}/E_{{7}}\cdot
S^{1})$:

\begin{enumerate}
\item[(6.13)] $A^{\ast }(E_{{8}}/E_{{7}}\cdot S^{1})=\mathbb{Z}[y_{{1}},y_{{6%
}},y_{{10}},y_{{15}}]/\left\langle r_{{15}},r_{{20}},r_{{24}},r_{{30}},y_{{1}%
}g_{{29}}\right\rangle $,
\end{enumerate}

\noindent where $y_{{1}},y_{{6}},y_{{10}},y_{{15}}$ are the same Schubert
classes as those in Theorem 7, and where if we let $\pi _{{m}}:\mathbb{Z}[y_{%
{1}},y_{{6}},y_{{10}},y_{{15}}]^{(2m)}\rightarrow A^{2m}(E_{{8}}/E_{{7}%
}\cdot S^{1})$ be induced from $\{y_{{1}},y_{{6}},y_{{10}},y_{{15}}\}\subset
A^{\ast }(E_{{8}}/E_{{7}}\cdot S^{1})$, then (Lemma 8)

\begin{quote}
1) for $m=15,20,24,30$, the $r_{{m}}\in \ker \pi _{{m}}$ in (6.13) should
satisfy
\end{quote}

$\qquad r_{{15}}\mid _{{y}_{{1}}{=0}}$ $=2y_{{15}}$; $\quad r_{{20}}\mid _{{y%
}_{{1}}{=0}}$ $=3y_{{10}}^{2}$; $\ r_{{24}}\mid _{{y}_{{1}}{=0}}$ $=5y_{{6}%
}^{4}$;

$\qquad r_{{30}}\mid _{{y}_{{1}}{=0}}$ $=y_{{15}}^{2}+y_{{10}}^{3}+y_{{6}%
}^{5}$

\begin{quote}
2) $\pi (g_{{29}})=s_{{29,1}}-s_{{29,2}}-s_{{29,3}}+\ s_{{29,4}}-s_{{29,5}%
}+s_{{29,6}}-s_{{29,7}}+s_{{29,8}}$.$\quad $
\end{quote}

With respect to the ordered basis $B(2m)$ of $\mathbb{Z}[y_{{1}},y_{{6}},y_{{%
10}},y_{{15}}]^{(2m)}$, $m=15,20,24,30$, the structure matrices $M(\pi _{{m}%
})$ have been computed by the \textsl{L--R coefficients} and their
corresponding Nullspaces $N(\pi _{{m}})$ are presented in [DZ$_{2}$, 7.5].
If we take, in view of Lemma 3, that

\renewcommand{\myy}[2]{y {#1}}

\begin{quote}
$r_{{15}}=2y_{{15}}-10y_{{1}}^{3}y_{{6}}^{2}\ -16y_{{1}}^{5}y_{{10}}+10y_{{1}%
}^{9}y_{{6}}\ -y_{{1}}^{15}$($=-\kappa _{{1}}$ in $N(\pi _{{15}})$);

$r_{{20}}=3y_{{10}}^{2}+10y_{{1}}^{2}y_{{6}}^{3}+18y_{{1}}^{4}y_{{6}}y_{{10}%
}-2y_{{1}}^{5}y_{{15}}-8y_{{1}}^{8}y_{{6}}^{2}+4y_{{1}}^{10}y_{{10}}-y_{{1}%
}^{14}y_{{6}}$

\qquad ($=-\kappa _{{2}}$ in $N(\pi _{{20}})$);

$r_{{24}}=5y_{{6}}^{4}+30y_{{1}}^{2}y_{{6}}^{2}y_{{10}}+15y_{{1}}^{4}y_{{10}%
}^{2}-2y_{{1}}^{9}y_{{15}}-5y_{{1}}^{12}y_{{6}}^{2}+y_{{1}}^{14}y_{{10}}$

\qquad ($=\frac{1}{2}\kappa _{{3}}-\frac{5}{2}\kappa _{{4}}$ in $N(\pi _{{24}%
})$);

$r_{{30}}=y_{_{{15}}}^{2}-8y_{_{10}}^{3}+y_{_{{6}}}^{5}-2y_{_{{1}}}^{3}y_{_{{%
6}}}^{2}y_{_{{15}}}+3y_{_{{1}}}^{4}y_{{6}}y_{{10}}^{2}-8y_{{1}}^{5}y_{{10}%
}y_{{15}}+6y_{{1}}^{9}y_{{6}}y_{{15}}$

$-9y_{{1}}^{10}y_{{10}}^{2}-y_{{1}}^{12}y_{{6}}^{3}-2y_{{1}}^{14}y_{{6}}y_{{%
10}}-3y_{{1}}^{15}y_{{15}}+8y_{{1}}^{20}y_{{10}}+y_{{1}}^{24}y_{{6}}-y_{{1}%
}^{30}$

\qquad ($=-\frac{1}{988483}\kappa _{{1}}+\frac{1}{988483}\kappa _{{2}}+\frac{%
8}{988483}\kappa _{{3}}-\frac{3}{988483}\kappa _{{5}}-\frac{2}{988483}\kappa
_{{6}}$

$\qquad -\frac{1}{988483}\kappa _{{7}}-\frac{9}{988483}\kappa _{{8}}+\frac{3%
}{988483}\kappa _{{9}}$ in $N(\pi _{{30}})$),
\end{quote}

\noindent then condition 1) is met by the set $\{r_{{15}},r_{{20}},r_{{24}%
},r_{{30}}\}$ of polynomials.

The proof will be completed once we show

\begin{enumerate}
\item[(6.14)] $y_{{1}}g_{{29}}\in \left\langle r_{{15}},r_{{20}},r_{{24}},r_{%
{30}}\right\rangle $.
\end{enumerate}

\noindent For this purpose we examine, in view of (6.13), the quotient map (%
\S 3.3)

\begin{quote}
$\varphi :\mathbb{Z}[y_{{1}},y_{{6}},y_{{10}},y_{{15}}]/\left\langle r_{{15}%
},r_{{20}},r_{{24}},r_{{30}}\right\rangle \rightarrow A^{\ast }(E_{{8}}/E_{{7%
}}\cdot S^{1})$.
\end{quote}

\noindent With $r_{{15}},r_{{20}},r_{{24}},r_{{30}}$ being obtained
explicitly, one finds that

\begin{quote}
$b(60)=18$; $\delta _{{60}}(r_{{15}},r_{{20}},r_{{24}},r_{{30}})=11$
(Example 4).
\end{quote}

\noindent On the other hand, granted with the basis theorem, we read from [DZ%
$_{{2}}$, 7.1] that rank($A^{60}$)$=7$. (6.14) is verified by Lemma 4.$%
\square $

\section{Remarks}

\textbf{7.1.} \textbf{Spectral sequence method (classical approaches).} In
order to compute the cohomology of a homogeneous space $G/H$, A. Borel, J.
Leray and H. Toda used the Leray--Serre spectral sequence for the fibration

\begin{center}
$G\hookrightarrow G/H\overset{\pi }{\rightarrow }BH$,
\end{center}

\noindent whereas P. Baum, D. Husemoller et al, and J. Wolf used the
Eilenberg--Moore spectral sequence for the fibration

\begin{center}
$G/H\overset{\pi }{\hookrightarrow }BH\overset{Bi}{\rightarrow }BG$,
\end{center}

\noindent where $BG$ is the classifying space of $G$, and where $Bi$ is the
induced map of the inclusion $i:H\rightarrow G$. When applying to the
integral cohomology of $G/H$ the computations start with the $E_{2}$--pages
([M, p.133; p.232])

\begin{center}
$E_{2}^{\ast ,\ast }=H^{\ast }(BH;H^{\ast }(G))$ or $E_{2}^{\ast ,\ast
}=Tor_{H^{\ast }(BG)}^{\ast ,\ast }(\mathbb{Z},H^{\ast }(BH))$,
\end{center}

\noindent respectively, where the integral cohomologies of the spaces $%
G,BG,BH$ are explicitly required as the inputs. We note that \textsl{the
integral cohomologies of the exceptional Lie groups, and of their
classifying spaces, have not yet been determined completely.}

On the other hand, our approach reduces the computation to the \textsl{%
Cartan numbers} of $G$ without resorting to either $H^{\ast }(G)$ or $%
H^{\ast }(BG)$.

\textbf{7.2. Generalizations.} The Grassmannians are special cases of flag
varieties. This raises the question whether our approach remains effective
for flag varieties of general types. In the subsequent work [DZ$_{3}$] the
method and results of this paper were applied to obtain the integral
cohomology rings of all complete flag manifolds associated to the
exceptional groups, in which the Chow rings of Grassmannians played an
important role.

\textbf{7.3. Computational aspects}. In appearance, the basis theorem (Lemma
2) and the \textsl{L--R coefficients} along may solve Problems 1 and 2 in \S %
3: one may start with many Schubert classes sufficient to generate the ring $%
A^{\ast }(G/H)$; specify the corresponding relations in \textsl{every degree}
using Nullspace (Lemma 3); and eliminate the extra generators (resp.
relations) using certain computer algebra packages. However, in practice,
computing with \textsl{L--R coefficients} tends to be extremely time
consuming as the degree $r$ running higher. This phenomenon is further
exaggerated by the fact that the rank of $A^{r}(G/H)$ (i.e. the number of
Schubert classes in degree $r$) turns to be very large as $r$ increasing to
the middle dimension $\frac{1}{2}\dim _{\mathbb{C}}G/H$. Similar difficulty
has been encountered by many authors, see discussion in Nikolenko and
Semenov [NS, \S 3].

Lemmas 7 and 8 are useful in reducing the computation cost considerably.
They specify \textsl{a minimal set }of generators and\textsl{\ the initial
constraints }(Lemma 8) of the corresponding\textsl{\ }relations on\ $A^{\ast
}(G/H)$ from that of the much simpler ring $H^{\ast }(G/H_{s})$ (see the
tables in \S 5). As results, computing with \textsl{L--R coefficients} was
limited to the relevant degrees only, and can be carried out on a personal
computer.

\bigskip

In problems 1 and 2, the minimum request on the generators and relations
comes from further geometric concerns. For instance, since $G/H$ is the
classifying space for all principal $H$--bundles whose $G$ reduction are
trivial, a minimal set of generators for $H^{\ast }(G/H)$ constitutes a set
of independent characteristic classes for those bundles. A minimal set of
relations with respect to the generators is also valuable in constructing
independent secondary characteristic classes [D$_{2}$].

\textbf{7.4.} \textbf{Relevant works.} In [Co, 1964] Conlon computed the
ring $H^{\ast }(E_{{6}}/D_{{5}})$ and the additive homology of $E_{{6}}/D_{{5%
}}\cdot S^{1}$. His method amounts to applying Morse theory to the space $%
\Omega (E_{{6}}/D_{{5}}\cdot S^{1},x,W)$ of paths to yield a cell
decomposition of $E_{{6}}/D_{{5}}\cdot S^{1}$ relative to $W$ in dimensions
less than $32$, here $W$ is the Cayley projective plane canonically embedded
in $E_{{6}}/D_{{5}}\cdot S^{1}$. Indeed, the basis theorem (Lemma 2) implies
the additive homology of any flag variety $G/H$.

In [IM, 2005] Iliev and Manivel described $A^{\ast }(E_{{6}}/D_{{5}}\cdot
S^{1})$ in terms of three Schubert classes and three relations ([IM,
Proposition 5.1--5.2]) by using divided difference operators due to Demazure
and Bernstein-Gelfand-Gelfand [D; BGG]. Our Theorem 4 indicates that two
Schubert classes and two relations suffice to present the ring. Early in
1974, Toda and Watanabe [TW, Corollary C] presented the ring $H^{\ast }(E_{{6%
}}/D_{{5}}\cdot S^{1})$ by two generators, although they did not describe it
in terms of Schubert classes.

In [NS, 2006], Nikolenko and Semenov investigated $A^{\ast }(E_{{8}}/E_{{7}%
}\cdot S^{1})\otimes \mathbb{Q}$. They specified generators in dimension $%
1,6 $ and $10$ in terms of Schubert classes, while the relations were given
by a number of equations that express monomials in the generators as linear
combinations of Schubert classes. Indeed, the ring $A^{\ast }(E_{{8}}/E_{{7}%
}\cdot S^{1})$ is more attractive than its rational analogue: $y_{{15}}$ in
Theorem 7 survives to the generator of the Chow ring of $E_{{8}}$ in degree $%
15$, which can not be detected from $A^{\ast }(E_{{8}}/E_{{7}}\cdot
S^{1})\otimes \mathbb{Q}$. In addition, even with rational coefficients,
presenting the relations by a minimal set of polynomials has further
implications. For instance, it is straightforward from Theorem 7 that, if we
let $g_{{20}},$ $g_{{24}},g_{{30}}$ be the polynomials obtained from the $r_{%
{20}},r_{{24}},r_{{30}}$ by substituting the ${y}_{{15}}$ with $8y_{{1}%
}^{5}y_{{10}}+5y_{{1}}^{3}y_{{6}}^{2}-5y_{{1}}^{9}y_{{6}}+\frac{1}{2}y_{{1}%
}^{15}$ by the relation $r_{15}$, then

\begin{quote}
$A^{\ast }(E_{{8}}/E_{{7}}\cdot S^{1})\otimes \mathbb{Q=Q}[y_{{1}},y_{{6}%
},y_{{10}}]/\left\langle g_{{20}},g_{{24}},g_{{30}}\right\rangle $.
\end{quote}

\noindent Consequently, the rational homotopy groups of $E_{{8}}/E_{{7}%
}\cdot S^{1}$ are given by

\begin{quote}
$\pi _{{r}}(E_{{8}}/E_{{7}}\cdot S^{1})\otimes \mathbb{Q}=\left\{ 
\begin{tabular}{l}
$\mathbb{Q}\text{ for }r=2,12,20,39,47,59$ \\ 
$0\text{ otherwise,}$%
\end{tabular}%
\right. $
\end{quote}

\noindent where the generators for the nontrivial $\pi _{{r}}\otimes \mathbb{%
Q}$ are fashioned from $y_{{1}},y_{{6}},y_{{10}},g_{{20}},$ $g_{{24}},g_{{30}%
}$ respectively ([BT, p.258-265]). In general, letting $G/H$ be one of the
Grassmannians concerned in this paper, one deduces directly from Theorems
1--7, using the method illustrated in [BT, p.258-265], that

\bigskip

\noindent \textbf{Corollary.} One has either i) $\pi _{{r}}(G/H)\otimes 
\mathbb{Q}=0$ or ii) $\pi _{{r}}(G/H)\otimes \mathbb{Q\cong Q}$, where ii)
occurs if and only if $r$ takes the values in the table below

\begin{center}
\begin{tabular}{l|l}
\hline\hline
$G/H$ & $\text{the }r\text{ with }\pi _{{r}}(G/H)\otimes \mathbb{Q\cong Q}$
\\ \hline
$F_{{4}}/C_{{3}}\cdot S^{1}$ & $2,8,15,23$ \\ \hline
$F_{{4}}/B_{{3}}\cdot S^{1}$ & $2,8,15,23$ \\ \hline
$E_{{6}}/A_{{6}}\cdot S^{1}$ & $2,6,8,15,17,23$ \\ \hline
$E_{{6}}/D_{{5}}\cdot S^{1}$ & $2,8,17,23$ \\ \hline
$E_{{7}}/E_{{6}}\cdot S^{1}$ & $2,10,18,19,27,35$ \\ \hline
$E_{{7}}/D_{{6}}\cdot S^{1}$ & $2,8,12,23,27,35$ \\ \hline
$E_{{8}}/E_{{7}}\cdot S^{1}$ & $2,12,20,39,47,59$ \\ \hline\hline
\end{tabular}
\end{center}

Historically, the mod $p$ cohomologies $H^{\ast }(G;\mathbb{Z}_{{p}})$ of
exceptional $G$ were achieved using case by case calculations (see Ka\v{c}
[K] for a thorough summary on the history). In 1974, H. Toda [T] initiated
the project computing the integral cohomology of homogeneous spaces $G/H$
with $G$ an exceptional Lie group and $H\subset G$ a torsion free subgroup
of maximal rank. This amounts to combining Borel's method [B$_{{1}}$] with
the previous results on $H^{\ast }(G;\mathbb{Z}_{{p}})$ (as a module over
the Steenrod algebra) for all primes $p$. After Toda, the cohomologies of
the $G/H$ considered in Theorems 1, 3--6 have been studied by Toda,
Watanabe, Ishitoya and Nakagawa in [I, IT, TW, W$_{{1}}$, W$_{{2}}$, N] in
which the generators are specified only up to degrees.

In comparison our approach is free of types. By taking generators among
Schubert classes on $G/H$ their geometric configurations are transparent in
view of the unified construction (2.4) of all Schubert varieties. Moreover,
instead of resorting to the $H^{\ast }(G;\mathbb{Z}_{{p}})$, this work
brings a way to determine $H^{\ast }(G;\mathbb{Z})$ [DZ$_{4}$].

\textbf{7.5. Corrections.} Theorem 3 corrects a mistake occurring in [W,
1998]. Toda and Ishitoya claimed in [IT, 1977] that the ring $H^{\ast }(E_{{6%
}}/A_{{6}}\cdot S^{1})$ is the quotient of a polynomial ring in eight
variables modulo an ideal generated by eight polynomials (with those eight
polynomials not being computed explicitly). Watanabe asserted in [W, 1998]
that it was generated by three elements in degrees $2,6$ and $8$
respectively. However, according to the proof of Theorem 3, \textsl{four} is
the minimal number of generators for this ring.

After the work [IT] Ishitoya made explicit computation about $H^{\ast }(E_{{6%
}}/A_{{6}}\cdot S^{1})$ and $H^{\ast }(E_{{6}}/A_{{6}})$ in [I] which
contains also an error. Corollary 3.5 in [I] implies that $H^{22}(E_{{6}}/A_{%
{6}})=0$. However, from the table in the proof of Theorem 10 one reads that $%
H^{22}(E_{{6}}/A_{{6}})=\mathbb{Z}_{3}$. These issues witness the subtleness
in the traditional approach.

\bigskip

\noindent \textbf{Acknowledgement.} The authors are grateful to their
referees for valuable suggestions and many improvements on the earlier
version of this paper.

Thanks are also due to Mamoru Mimura, Nakagawa, Totaro and Zainoulline for
communications concerning the work. In particular, Nakagawa recently
announced a presentation for the ring $H^{\ast }(E_{{8}}/E_{{7}}\cdot S^{1})$%
.

\begin{center}
\textbf{References}
\end{center}

\begin{enumerate}
\item[{[B]}] P. Baum, On the cohomology of homogeneous spaces, Topology
7(1968), 15-38.

\item[{[B$_{{1}}$]}] A. Borel, Sur la cohomologie des espaces fibr\'{e}s
principaux et des espaces homogenes de groupes de Lie compacts, Ann. Math.
57(1953), 115-207.

\item[{[BGG]}] I. N. Bernstein, I. M. Gelfand, S. I. Gelfand, Schubert cells
and cohomology of the spaces G/P, Russian Math. Surveys 28 (1973), 1-26.

\item[{[BH]}] A. Borel and F. Hirzebruch, Characteristic classes and
homogeneous spaces (I), Amer. J. Math. 80, 1958, 458--538.

\item[{[BS]}] R. Bott and H. Samelson, Application of the theory of Morse to
symmetric spaces, Amer. J. Math., Vol. LXXX, no. 4 (1958), 964-1029.

\item[{[BT]}] R. Bott and L. Tu, Differential forms in algebraic topology.
Graduate Texts in Mathematics, 82. Springer-Verlag, New York-Berlin, 1982.

\item[{[C]}] C. Chevalley, Sur les d\'{e}compositions celluaires des Espaces
G/B, in Algebraic groups and their generalizations: Classical methods, W.
Haboush ed. Proc. Symp. in Pure Math. 56 (part 1) (1994), 1-26.

\item[{[Co]}] L. Conlon, On the topology of $EIII$ and $EIV$, Proc. Amer.
Math. Soc., 16(1965), 575-581.

\item[{[D$_{1}$]}] H. Duan, Multiplicative rule of Schubert classes, Invent.
Math.159 (2005), 407-436.

\item[{[D$_{2}$]}] H. Duan, Characteristic classes for complex bundles whose
real reductions are trivial, Proc. Amer. Math. Soc. Vol. 128 (2000),
2465-2471.

\item[{[DZ$_{{1}}$]}] H. Duan and X. Zhao, Algorithm for multiplying Schubert
classes. Internat. J. Algebra Comput. 16(2006), 1197--1210.

\item[{[DZ$_{{2}}$]}] H. Duan and X. Zhao, Appendix to \textquotedblleft The
Chow rings of generalized Grassmannians\textquotedblright , arXiv:
math.AG/0510085.

\item[{[DZ$_{{3}}$]}] H. Duan and X. Zhao, The integral cohomology of
complete flag manifolds, arXiv: math.AT/0801.2444

\item[{[DZ$_{{4}}$]}] H. Duan and X. Zhao, The cohomology of Lie groups,
arXiv: math.AT /0711.2541.

\item[{[Fu]}] W. Fulton, Intersection theory, Springer--Verlag, 1998.

\item[{[H]}] H.C. Hansen, On cycles in flag manifolds, Math. Scand. 33
(1973), 269-274.

\item[{[Hi]}] H. Hiller, Geometry of Coxeter groups, \textsl{Research Notes
in Mathematics}, 54. Pitman Advanced Publishing Program, 1982.

\item[{[HMS]}] D. Husemoller, J. Moore, J. Stasheff, Differential homological
algebra and homogeneous spaces, J. Pure Appl. Algebra 5 (1974), 113--185.

\item[{[Hu]}] J. E. Humphreys, Introduction to Lie algebras and
representation theory, Graduated Texts in Math. 9, Springer-Verlag New York,
1972.

\item[{[IM]}] A. Iliev and L. Manivel, The Chow ring of the Cayley plane,
Compositio Math., 141(2005), 146-160.

\item[{[I]}] K. Ishitoya, Integral cohomology ring of the symmetric space
EII. J. Math. Kyoto Univ. 17(1977), no.2, 375--397.

\item[{[IT]}] K. Ishitoya, H. Toda, On the cohomology of irreducible
symmetric spaces of exceptional type. J. Math. Kyoto Univ. 17(1977),
225--243.

\item[{[K]}] V.G. Kac, Torsion in cohomology of compact Lie groups and Chow
rings of reductive algebraic groups, Invent. Math. 80(1985), 69-79.

\item[{[M]}] J. McCleary, A user's guide to spectral sequences, Cambridge
University Press, 2001.

\item[{[MS]}] J. Milnor and J. Stasheff, Characteristic classes, Ann. of
Math. Studies 76, Princeton Univ. Press, 1975.

\item[{[N]}] M. Nakagawa, The integral cohomology ring of $E_{{7}}/T$. J.
Math. Kyoto Univ. 41(2001), 303--321.

\item[{[NS]}] S. I. Nikolenko and N. S. Semenov, Chow ring structure made
simple, arXiv: math.AG/0606335.

\item[{[S]}] J. Scherk, Algebra. A computational introduction, \textsl{%
Studies in Advanced Mathematics}, Chapman \& Hall/CRC, Boca Raton, FL, 2000.

\item[{[T]}] H. Toda, On the cohomology ring of some homogeneous spaces. J.
Math. Kyoto Univ. 15(1975), 185--199.

\item[{[TW]}] H. Toda, T. Watanabe, The integral cohomology ring of $F_{{4}%
}/T $ and $E_{{6}}/T$, J. Math. Kyoto Univ. 14(1974), 257--286.

\item[{[W$_{{1}}$]}] T. Watanabe, The integral cohomology ring of the
symmetric space EVII. J. Math. Kyoto Univ. 15(1975), 363--385

\item[{[W$_{{2}}$]}] T. Watanabe, Cohomology of the homogeneous space $%
E_{6}/T^{1}\cdot SU(6)$, \textsl{Group representations: cohomology, group
actions and topology}, 511--518, Proc. Sympos. Pure Math., 63, Amer. Math.
Soc., Providence, RI, 1998.

\item[{[Wo]}] J. Wolf, The cohomology of homogeneous spaces. Amer. J. Math.
99 (1977), no. 2, 312--340.
\end{enumerate}

\end{document}